\documentclass[11pt]{amsart}
\usepackage{amsmath,amssymb,amsthm,amsfonts,mathrsfs}
\usepackage{needspace}
\usepackage{enumitem}
\usepackage[margin=1in]{geometry}
\usepackage[english]{babel}
\usepackage[colorlinks,linkcolor=blue,citecolor=blue,urlcolor=blue]{hyperref}
\numberwithin{equation}{section}
\newtheorem{theorem}{Theorem}[section]
\newtheorem{proposition}[theorem]{Proposition}
\newtheorem{lemma}[theorem]{Lemma}
\newtheorem{corollary}[theorem]{Corollary}
\theoremstyle{definition}

\newcommand{\R}{\mathbb R}
\newcommand{\Z}{\mathbb Z}

\newcommand{\tr}{\operatorname{tr}}
\newcommand{\Id}{\operatorname{Id}}
\newcommand{\Hom}{\operatorname{Hom}}
\newcommand{\Hess}{\operatorname{Hess}}
\newcommand{\ii}{i}
\DeclareMathOperator{\Ric}{Ric}
\title[Sharp fill-in inequalities for spin initial data sets]{Sharp fill-in inequalities and rigidity for spin initial data sets}
\author{Tsz-Kiu Aaron Chow}
\address{Department of Mathematics, The Hong Kong University of Science and Technology, Clear Water Bay, Kowloon, Hong Kong}
\email{chowtka@ust.hk}
\thanks{The author is partially supported by the Hong Kong General Research Fund \#26310426 of the Hong Kong Research Grants Council and the National Natural Science Foundation of China under the Young Scientists Fund (Category C), Grant No.~NSFC26SC44.}
\date{}

\begin{document}

\begin{abstract}
We prove sharp fill-in inequalities for compact spin initial data sets satisfying the dominant energy condition. The inequalities relate the spacetime mean curvature of the boundary to its hyperspherical radius and stable one-systole. Equality characterizes spacelike graphs of constant boundary height over round Euclidean balls and, in the systolic case after passing to the universal cover, over products of a round disk with Euclidean space. We also prove rigidity under comparison with strictly convex Euclidean domains, using normalized boundary metrics. The proofs use a modified Dirac operator with a chiral boundary condition.
\end{abstract}

\maketitle

\section{Introduction}

Fill-in problems provide a local approach to global questions in scalar curvature geometry, connecting the boundary data of compact regions with results such as the positive mass theorem on noncompact manifolds. Scalar curvature bounds on a filling lead to sharp inequalities relating the mean curvature of its boundary to the hyperspherical radius and boundary systoles \cite{BrendleHung,BrendleTsiamisWang, GromovLectures}. In this paper, we prove corresponding fill-in inequalities for spin initial data sets satisfying the dominant energy condition, with the spacetime mean curvature of the boundary replacing its Riemannian mean curvature. We also characterize equality.

For a closed connected oriented Riemannian manifold $\Sigma^{n-1}$, its hyperspherical radius $\operatorname{Rad}(\Sigma)$ is the supremum of all $r>0$ for which there is a Lipschitz map $N:\Sigma\to S^{n-1}$ with nonzero degree and $\operatorname{Lip}(N)\leq r^{-1}$. Let $(M^n,g)$ be a compact connected spin manifold with boundary $\partial M=\Sigma$, and let $H_\Sigma$ denote the mean curvature, with convention that Euclidean spheres having positive mean curvature. Gromov's fill-in inequality \cite[Section~4.3]{GromovLectures} asserts that $R_g\geq0$ implies
\[
\min_\Sigma H_\Sigma\leq(n-1)/\operatorname{Rad}(\Sigma).
\]
Equality characterizes Euclidean balls \cite[Theorem~1.5]{CHZ}. For $R_g\geq-n(n-1)$, Brendle, Tsiamis, and Wang \cite[Corollary~2]{BrendleTsiamisWang} proved the corresponding bound $\min_\Sigma H_\Sigma\leq(n-1)\sqrt{1+\operatorname{Rad}(\Sigma)^{-2}}$ in dimensions $n\geq3$. The equality case characterizes hyperbolic geodesic balls \cite[Theorem~A]{AmmannLockman}.

Both scalar curvature assumptions above arise naturally from the dominant energy condition. An initial data set consists of a Riemannian manifold $(M^n,g)$ and a symmetric two-tensor $k$, with energy and momentum densities
\[
\mu=\frac12\bigl(R_g+(\tr_g k)^2-|k|_g^2\bigr),\qquad J=\operatorname{div}_g k-d\tr_g k.
\]
The dominant energy condition is $\mu\geq|J|_g$. It reduces to $R_g\geq0$ when $k=0$ and to $R_g\geq-n(n-1)$ when $k=g$. Our first result generalizes Gromov's fill-in inequality to initial data sets.

\begin{theorem}\label{thm:radius}
Let $(M^n,g,k)$, $n\geq2$, be a compact connected spin initial data set satisfying $\mu\geq|J|_g$, with connected boundary $\partial M=\Sigma$. Then
\[
\min_\Sigma\left(H_\Sigma-\sqrt{(\tr_\Sigma k)^2+\frac{(n-1)^2}{\operatorname{Rad}(\Sigma)^2}}\right)\leq0.
\]
Equality holds if and only if there is $f\in C^\infty(M)$, constant on $\Sigma$, such that $k=(1+|df|_g^2)^{-1/2}\Hess_g f$ and $(M,g+df\otimes df)$ is isometric to the closed Euclidean ball of radius $\operatorname{Rad}(\Sigma)$.

Equivalently, equality holds if and only if $(M,g,k)$ is a spacelike graph of constant boundary height over a round Euclidean ball in Minkowski space.
\end{theorem}

The choices $k=0$ and $k=g$ recover the Euclidean and hyperbolic inequalities and rigidity statements above. Under $H_\Sigma\geq|\tr_\Sigma k|$, the estimate in Theorem~\ref{thm:radius} is equivalent to 
\[\min_\Sigma\sqrt{H_\Sigma^2-(\tr_\Sigma k)^2}\leq(n-1)/\operatorname{Rad}(\Sigma),\]
where $\sqrt{H_\Sigma^2-(\tr_\Sigma k)^2}$ is the norm of the spacetime mean curvature vector. Raulot \cite[Corollary~2]{RaulotFillin} proved this inequality for $3\leq n\leq7$ under $H_\Sigma>|\tr_\Sigma k|$. His proof uses a complete Jang graph and a Dirac estimate with MIT boundary conditions. Our proof works directly on the initial data set, using a modified Dirac operator with a chiral boundary condition. It applies in every dimension $n\geq2$. For related spectral estimates, see \cite{BarRadius,HMR,RaulotDEC}.

Boundary mean curvature also controls the length of nontrivial boundary curves. Under a topological assumption modelled on $B^2\times\mathbb T^{n-2}$, Brendle and Hung \cite[Theorem~1.4]{BrendleHung} proved that a compact Riemannian manifold of dimension $3\leq n\leq7$ with $R_g\geq0$ and $\min_\Sigma H_\Sigma>0$ has a homologically nontrivial closed boundary curve of length at most $2\pi/\min_\Sigma H_\Sigma$. Equality in their systolic inequality forces the universal cover to be a flat disk product \cite[Theorem~1.1]{Wang}. The proofs use minimal slicing. For earlier results in dimension three, see \cite{Ambrozio,Longa}.

We prove a spacetime analogue for the stable one-systole of the boundary. For $\alpha\in H_1(\Sigma;\R)$, the stable norm $\|\alpha\|$ is the infimum of $\sum_j|r_j|\operatorname{Length}(\gamma_j)$ over real Lipschitz one-cycles $\sum_jr_j\gamma_j$ representing $\alpha$. The stable one-systole introduced in  \cite{BangertKatz} is defined as
\[
\operatorname{stsys}_1(\Sigma)=\min\{\|\alpha\|:0\neq\alpha\in\operatorname{im}(H_1(\Sigma;\Z)\to H_1(\Sigma;\R))\}.
\]
Write $B^2(r)$ for the closed Euclidean disk of radius $r$.

\begin{theorem}\label{thm:stable-systole}
Let $(M^n,g,k)$, $n\geq3$, be a compact connected spin initial data set satisfying $\mu\geq|J|_g$, with connected boundary $\partial M=\Sigma$ and $b_1(\Sigma)=n-1$. Suppose that there are smooth maps $\Theta:M\to\mathbb T^{n-2}$ and $\theta_0:\Sigma\to S^1$ such that $\deg(\theta_0,\Theta|_\Sigma)\neq0$. Then
\[
\min_\Sigma\left(H_\Sigma-\sqrt{(\tr_\Sigma k)^2+\frac{4\pi^2}{\operatorname{stsys}_1(\Sigma)^2}}\right)\leq0.
\]
Equality holds if and only if there is $f\in C^\infty(M)$, constant on $\Sigma$, such that $k=(1+|df|_g^2)^{-1/2}\Hess_g f$ and the universal cover of $(M,g+df\otimes df)$ is isometric to
\[
B^2\!\left((2\pi)^{-1}\operatorname{stsys}_1(\Sigma)\right)\times\R^{n-2}
\]
with the standard Euclidean product metric.

Equivalently, equality holds if and only if the universal cover of $M$, with the lifted initial data, is a spacelike graph of constant boundary height over this disk product in Minkowski space.
\end{theorem}

The flat products $B^2(r)\times S^1(L)^{n-2}$, with $k=0$ and $L\geq r$, attain equality. Here $S^1(L)$ is the circle of radius $L$, so $H_\Sigma=1/r$ and $\operatorname{stsys}_1(\Sigma)=2\pi r$. When $k=0$, the function $f$ in the equality statement is constant. In dimensions three through seven, the Riemannian inequality and rigidity also follow from \cite{BrendleHung,Wang}, since their shortest curve has length at least $\operatorname{stsys}_1(\Sigma)$. We retain their topological assumption and require in addition that $b_1(\Sigma)=n-1$, while treating general initial data in every dimension $n\geq3$ under the spin assumption.

The main ingredient in the proofs of Theorems~\ref{thm:radius} and~\ref{thm:stable-systole} is the following generalization of a rigidity theorem proved in joint work of the author with B\"ar, Brendle, and Hanke \cite[Theorem~I]{BBCH}. Write $\|dN\|_{\tr}$ for the sum of the singular values of $dN$.

\begin{theorem}\label{thm:boundary-rigidity}
Let $(M^n,g,k)$, $n\geq2$, be a compact connected spin initial data set with nonempty boundary satisfying $\mu\geq|J|_g$. Suppose that a smooth map $N:\Sigma\to S^{n-1}$ has positive degree and
\[
H_\Sigma\geq\sqrt{(\tr_\Sigma k)^2+\|dN\|_{\tr}^2}.
\]
Then there is $f\in C^\infty(M)$, constant on $\Sigma$, such that $k=(1+|df|_g^2)^{-1/2}\Hess_g f$ and $(M,g+df\otimes df)$ is isometric to a compact convex Euclidean domain with smooth boundary. Under this isometry, $N$ agrees with the outward Gauss map up to an orthogonal transformation, and $H_\Sigma=\sqrt{(\tr_\Sigma k)^2+\|dN\|_{\tr}^2}$ everywhere. Conversely, every such initial data set with its outward Gauss map satisfies these hypotheses.

Equivalently, $(M,g,k)$ is a spacelike graph of constant boundary height over a compact convex Euclidean domain in Minkowski space.
\end{theorem}

We note that the boundary assumption in \cite[Theorem~I]{BBCH} is $H_\Sigma\geq|\tr_\Sigma k|+\|dN\|_{\tr}$. The height function of the resulting graph is constant on the boundary and has zero normal derivative wherever the Euclidean second fundamental form is nonzero \cite[Example~1.2]{BBCH}. On the other hand, Theorem~\ref{thm:boundary-rigidity} characterize all smooth spacelike graphs of constant boundary height over compact convex Euclidean domains. Indeed, the graph and its domain $K$ induce the same boundary metric. If $N_0$ and $H_0$ are the Euclidean Gauss map and mean curvature of $\partial K$, then
\[
H_\Sigma^2-(\tr_\Sigma k)^2=H_0^2=\|dN_0\|_{\tr}^2.
\]
Related boundary index and rigidity results for nonnegative scalar curvature appear in \cite{LottBoundary,WangXie,CHZ}.

Theorem~\ref{thm:boundary-rigidity} also yields a spacetime analogue of Llarull's theorem for compact initial data sets with boundary. For closed manifolds, Llarull's theorem \cite{Llarull} states that the unit round metric $g_{\mathrm{round}}$ on $S^n$ is the only metric $g$ satisfying $g\geq g_{\mathrm{round}}$ and $R_g\geq n(n-1)$. Min-Oo's conjecture, disproved in \cite{BMN}, prescribed the same scalar curvature lower bound and a unit round totally geodesic boundary, without requiring interior metric comparison. With the metric, scalar curvature, and boundary mean curvature bounded below by their round model values, rigidity holds for three-dimensional spherical caps \cite[Theorem~1.5]{HuLiuShi} and for hemispheres in all dimensions \cite[Corollary~1.1]{LeeTam}.

The following corollary gives a spacetime analogue for convex Euclidean domains, with the dominant energy condition in the interior and a metric inequality on the boundary after normalization by $H_\Sigma^2-(\tr_\Sigma k)^2$. For related Riemannian results using mean-curvature normalization, see \cite{YaoComparison,KoYao}.

\begin{corollary}\label{cor:normalized-comparison}
Let $(M^n,g,k)$, $n\geq3$, be a compact connected spin initial data set with connected boundary $\partial M=\Sigma$ satisfying $\mu\geq|J|_g$ and $H_\Sigma\geq0$. Let $K_0\subset\R^n$ be a compact convex domain with smooth boundary whose second fundamental form is positive definite, and let $H_0$ denote its mean curvature. Suppose that a smooth map $F:\Sigma\to\partial K_0$ has positive degree and
\[
\bigl(H_\Sigma^2-(\tr_\Sigma k)^2\bigr)g|_\Sigma\geq(H_0\circ F)^2F^*g_{\R^n}.
\]
Then there are $f\in C^\infty(M)$, constant on $\Sigma$, and $a>0$ such that $k=(1+|df|_g^2)^{-1/2}\Hess_g f$ and $F$ extends to an isometry from $(M,a^2(g+df\otimes df))$ onto $(K_0,g_{\R^n})$. In particular,
\[
\bigl(H_\Sigma^2-(\tr_\Sigma k)^2\bigr)g|_\Sigma=(H_0\circ F)^2F^*g_{\R^n}.
\]
Conversely, every such initial data set and orientation-preserving isometry satisfy these hypotheses. The factor $a$ is necessary because the normalized boundary metrics are scale invariant.
\end{corollary}

In Section~\ref{sec:boundary-estimate}, we prove the sharp estimate for the modified Dirac operator with the chiral boundary condition and combine it with the index argument in \cite{BBCH} to obtain parallel sections. In Section~\ref{sec:comparison-metric}, we use the spinor current to construct the metric $g+df\otimes df$, and we prove Theorem~\ref{thm:boundary-rigidity} and Corollary~\ref{cor:normalized-comparison}. Section~\ref{sec:radius} applies these results to sphere maps approaching the hyperspherical radius. In Section~\ref{sec:stable-systole}, stable norm duality gives sphere maps on finite covers and yields the systolic inequality. In the equality case, the limiting spinor current gives a Ricci-flat metric, and harmonic one-forms determine its flat disk-product universal cover.

\noindent\textit{Disclosure of AI tools.} The author developed the mathematical ideas and arguments presented in this paper. ChatGPT 6 Pro was used to assist with calculations, checking proof details, and revising the exposition. The author takes full responsibility for the mathematical content and the final manuscript.

\section{The boundary identity and the sharp estimate}\label{sec:boundary-estimate}

Let $(M^n,g,k)$ be a compact connected spin initial data set with $n\geq2$ and nonempty boundary $\partial M=\Sigma$. Let $\nu$ be the outward unit normal, and write $D$ for the Levi--Civita connection of $g$ and $h(X,Y)=g(D_X\nu,Y)$ for the second fundamental form of $\Sigma$, so that $H_\Sigma=\tr h$. Let $S_0$ be the spinor module of the complex Clifford algebra of $\R^n$, and let $\mathcal S$ be the spinor bundle over $M$. As in \cite{BBCH}, we define the complex vector bundle
\[
\mathcal E=\mathcal S\otimes S_0^*=\Hom(S_0,\mathcal S),
\]
with the spin connection on the first factor and the trivial connection on the second. Write $\gamma$ and $\gamma_0$ for the Clifford multiplication on $\mathcal S$ and $S_0$, and define the dual action by $\gamma_0^*(v)\ell=\ell\circ\gamma_0(v)$. Thus $(\gamma(X)\otimes\gamma_0^*(v))s=\gamma(X)\circ s\circ\gamma_0(v)$ for a section $s$ of $\mathcal E$.

Denote $k(X)=k(X,\cdot)^\sharp$. We define the connection on $\mathcal E$ by
\begin{equation}\label{eq:modified-connection}
\nabla_X^{\mathcal E,k}s=\nabla_X^{\mathcal E}s+\frac{\ii}{2}(\gamma(k(X))\otimes\Id)s.
\end{equation}
For a local orthonormal frame $e_1,\ldots,e_n$, its Dirac operator is
\[
\sum_{j=1}^n(\gamma(e_j)\otimes\Id)\nabla_{e_j}^{\mathcal E,k}=\mathcal D-\frac{\ii}{2}\tr_g k,
\]
where $\mathcal D$ denotes the Dirac operator of $\nabla^{\mathcal E}$.  For a smooth map $N:\Sigma\to S^{n-1}$, define $\chi: \mathcal E|_\Sigma\to \mathcal E|_\Sigma$ by
\[
\chi=-\gamma(\nu)\otimes\gamma_0^*(N).
\]
Then $\chi^*=\chi$ and $\chi^2=\Id$. Along $\Sigma$, the boundary Dirac operator $\mathcal D^\Sigma$ acting on sections of $\mathcal E|_\Sigma$ is defined by
\[
\mathcal D^\Sigma=\sum_{a=1}^{n-1}(\gamma(\nu)\gamma(e_a)\otimes\Id)\nabla^{\mathcal E}_{e_a}+\frac{H_\Sigma}{2},
\]
where $e_1,\ldots,e_{n-1}$ is a local orthonormal frame on $\Sigma$. Moreover, we adopt the convention that $a,b$ denote boundary indices ranging from $1$ to $n-1$.

\begin{lemma}\label{lem:boundary-identity}
For every $s\in H^1(M,\mathcal E)$ satisfying $\chi s=s$ on $\Sigma$, we have
\[
\begin{aligned}
&-\int_M\left|\left(\mathcal D-\frac{\ii}{2}\tr_g k\right)s\right|^2+\int_M|\nabla^{\mathcal E,k}s|^2+\frac12\int_M\langle(\mu+\ii\gamma(J)\otimes\Id)s,s\rangle\\
&=-\frac12\int_\Sigma H_\Sigma|s|^2-\frac12\int_\Sigma\left\langle\left(\sum_a\gamma(e_a)\otimes\gamma_0^*(dN(e_a))+\ii(\tr_\Sigma k)\gamma(\nu)\otimes\Id\right)s,s\right\rangle.
\end{aligned}
\]
\end{lemma}

\begin{proof}
We first assume that $s$ is smooth. By \cite[Corollary~3.4]{BBCH},
\[
\begin{aligned}
&-\int_M\left|\left(\mathcal D-\frac{\ii}{2}\tr_g k\right)s\right|^2+\int_M|\nabla^{\mathcal E,k}s|^2+\frac12\int_M\langle(\mu+\ii\gamma(J)\otimes\Id)s,s\rangle\\
&=\int_\Sigma\left\langle\left(\mathcal D^\Sigma-\frac{H_\Sigma}{2}+\frac{\ii}{2}\gamma(k(\nu)-(\tr_g k)\nu)\otimes\Id\right)s,s\right\rangle.
\end{aligned}
\]
Since $\chi s=s$ and $\chi^*=\chi$, \cite[Lemma~2.2]{BBCH} gives
\[
\begin{aligned}
2\langle\mathcal D^\Sigma s,s\rangle&=\operatorname{Re}\langle(\mathcal D^\Sigma\chi+\chi\mathcal D^\Sigma)s,s\rangle
=-\left\langle\sum_a(\gamma(e_a)\otimes\gamma_0^*(dN(e_a)))s,s\right\rangle.
\end{aligned}
\]
We now observe that for $X\in T\Sigma$, the Clifford relations give
\begin{equation}\label{eq:tangential-current}
\begin{aligned}
2\langle(\ii\gamma(X)\otimes\Id)s,s\rangle&=\left\langle\bigl((\ii\gamma(X)\otimes\Id)\chi+\chi(\ii\gamma(X)\otimes\Id)\bigr)s,s\right\rangle\\
&=-\ii\left\langle\bigl((\gamma(X)\gamma(\nu)+\gamma(\nu)\gamma(X))\otimes\gamma_0^*(N)\bigr)s,s\right\rangle\\
&=0.
\end{aligned}
\end{equation}
Using $k(\nu)-(\tr_g k)\nu=\sum_a k(\nu,e_a)e_a-(\tr_\Sigma k)\nu$, we obtain
\[
\begin{aligned}
&\langle(\ii\gamma(k(\nu)-(\tr_g k)\nu)\otimes\Id)s,s\rangle\\
&=\sum_a k(\nu,e_a)\langle(\ii\gamma(e_a)\otimes\Id)s,s\rangle-(\tr_\Sigma k)\langle(\ii\gamma(\nu)\otimes\Id)s,s\rangle\\
&=-(\tr_\Sigma k)\langle(\ii\gamma(\nu)\otimes\Id)s,s\rangle.
\end{aligned}
\]
Consequently,
\[
\begin{aligned}
&\left\langle\left(\mathcal D^\Sigma-\frac{H_\Sigma}{2}+\frac{\ii}{2}\gamma(k(\nu)-(\tr_g k)\nu)\otimes\Id\right)s,s\right\rangle\\
&=-\frac{H_\Sigma}{2}|s|^2-\frac12\left\langle\left(\sum_a\gamma(e_a)\otimes\gamma_0^*(dN(e_a))+\ii(\tr_\Sigma k)\gamma(\nu)\otimes\Id\right)s,s\right\rangle.
\end{aligned}
\]
Substituting into the first identity proves the assertion for smooth $s$.

Since the $+1$-eigenspaces of $\chi$ form a smooth subbundle, smooth sections satisfying $\chi s=s$ are dense in the $H^1$ space with the same boundary condition. The identity follows by continuity of the trace map $H^1(M,\mathcal E)\to L^2(\Sigma,\mathcal E|_\Sigma)$.
\end{proof}

\begin{lemma}\label{lem:sharp-boundary}
For every $s\in H^1(M,\mathcal E)$ satisfying $\chi s=s$ on $\Sigma$, we have
\[
\begin{aligned}
&-\int_M\left|\left(\mathcal D-\frac{\ii}{2}\tr_g k\right)s\right|^2+\int_M|\nabla^{\mathcal E,k}s|^2+\frac12\int_M\langle(\mu+\ii\gamma(J)\otimes\Id)s,s\rangle\\
&\leq\frac12\int_\Sigma\left(\sqrt{(\tr_\Sigma k)^2+\|dN\|_{\tr}^2}-H_\Sigma\right)|s|^2.
\end{aligned}
\]
Here $\|dN\|_{\tr}$ is the sum of the singular values of $dN$.
\end{lemma}

\begin{proof}
We first observe that
\begin{equation}\label{eq:boundary-anticommutator}
\begin{aligned}
&(\gamma(e_a)\otimes\gamma_0^*(dN(e_a)))(\ii\gamma(\nu)\otimes\Id)+(\ii\gamma(\nu)\otimes\Id)(\gamma(e_a)\otimes\gamma_0^*(dN(e_a)))\\
&=\ii(\gamma(e_a)\gamma(\nu)+\gamma(\nu)\gamma(e_a))\otimes\gamma_0^*(dN(e_a))&\\
&=0.
\end{aligned}
\end{equation}
Together with $(\ii\gamma(\nu)\otimes\Id)^2=\Id$, we get
\begin{equation}\label{eq:boundary-square}
\begin{aligned}
&\left(\sum_a\gamma(e_a)\otimes\gamma_0^*(dN(e_a))+\ii(\tr_\Sigma k)\gamma(\nu)\otimes\Id\right)^2\\
&=\left(\sum_a\gamma(e_a)\otimes\gamma_0^*(dN(e_a))\right)^2+(\tr_\Sigma k)^2\Id.
\end{aligned}
\end{equation}
\cite[Lemma~3.5]{BBCH} implies
\[
\begin{aligned}
\left|\sum_a(\gamma(e_a)\otimes\gamma_0^*(dN(e_a)))s\right|&\le\|dN\|_{\tr}|s|.
\end{aligned}
\]
Together with \eqref{eq:boundary-square} and the property that both endomorphisms $\gamma(e_a)\otimes\gamma_0^*(dN(e_a))$ and $\ii\gamma(\nu)\otimes\Id$ are self-adjoint, this gives
\begin{equation}\label{eq:boundary-operator-bound}
\begin{aligned}
&\left|\left(\sum_a\gamma(e_a)\otimes\gamma_0^*(dN(e_a))+\ii(\tr_\Sigma k)\gamma(\nu)\otimes\Id\right)s\right|^2\\
&=\left|\sum_a(\gamma(e_a)\otimes\gamma_0^*(dN(e_a)))s\right|^2+(\tr_\Sigma k)^2|s|^2\\
&\leq\bigl(\|dN\|_{\tr}^2+(\tr_\Sigma k)^2\bigr)|s|^2.
\end{aligned}
\end{equation}
Applying \eqref{eq:boundary-operator-bound} to Lemma~\ref{lem:boundary-identity}, we obtain
\begin{align*}
&-\int_M\left|\left(\mathcal D-\frac{\ii}{2}\tr_g k\right)s\right|^2+\int_M|\nabla^{\mathcal E,k}s|^2+\frac12\int_M\langle(\mu+\ii\gamma(J)\otimes\Id)s,s\rangle\\
&\leq-\frac12\int_\Sigma H_\Sigma|s|^2+\frac12\int_\Sigma\left|\left(\sum_a\gamma(e_a)\otimes\gamma_0^*(dN(e_a))+\ii(\tr_\Sigma k)\gamma(\nu)\otimes\Id\right)s\right||s|\\
&\leq\frac12\int_\Sigma\left(\sqrt{(\tr_\Sigma k)^2+\|dN\|_{\tr}^2}-H_\Sigma\right)|s|^2.
\end{align*}
\end{proof}

In the next step, we establish a general compactness statement for later applications in hyperspherical radius and systolic inequality.

\begin{lemma}\label{lem:compactness}
Let $\widetilde M\to M$ be a connected covering and let $\mathcal C\subset\widetilde M$ be a compact set with nonempty interior. Suppose that there are smooth sections $s_j$ of $\mathcal E$ on $\widetilde M$ satisfying
\[
\int_{\mathcal C}|s_j|^2=1,\qquad \int_K|\nabla^{\mathcal E,k}s_j|^2\longrightarrow0
\]
for every compact domain $K\subset\widetilde M$. Then a subsequence converges strongly in $H^1_{\mathrm{loc}}$, up to the boundary, to a smooth nowhere-zero section $s$ on $\widetilde M$ with $\nabla^{\mathcal E,k}s=0$. The boundary traces converge strongly in $L^2_{\mathrm{loc}}(\partial\widetilde M,\mathcal E)$.
\end{lemma}

\begin{proof}
	Fix a connected compact Lipschitz domain $K\subset\widetilde M$ containing $\mathcal C$. We apply \cite[Proposition~B.1]{BBCH} on finitely many coordinate Lipschitz domains covering $K$, starting with one contained in $\operatorname{int}\mathcal C$ and extending the estimate across their open overlaps. The arguments in \cite[proofs of Corollaries~6.4--6.5]{BBCH} gives
\[
\|v\|_{H^1(K)}^2+\|v\|_{L^2(K\cap\partial\widetilde M)}^2\leq C_K\left(\int_{\mathcal C}|v|^2+\int_K|\nabla^{\mathcal E,k}v|^2\right)
\]
for every smooth section $v$ on $K$.

This gives local $H^1$ bounds for $s_j$ up to the boundary. Rellich compactness and a diagonal subsequence give weak convergence in $H^1_{\mathrm{loc}}$ and strong convergence in $L^2_{\mathrm{loc}}$ to a section $s$ with $\int_{\mathcal C}|s|^2=1$ and $\nabla^{\mathcal E,k}s=0$. In particular, $s$ is smooth up to the boundary. If $s$ vanished at one point, uniqueness for the parallel equation along curves would imply $s\equiv0$, contrary to $\int_{\mathcal C}|s|^2=1$. Applying the estimate to $s_j-s$, as in \cite[proof of Proposition~6.8]{BBCH}, we obtain
\begin{align*}
\|s_j-s\|_{H^1(K)}^2+\|s_j-s\|_{L^2(K\cap\partial\widetilde M)}^2
&\leq C_K\left(\int_{\mathcal C}|s_j-s|^2+\int_K|\nabla^{\mathcal E,k}(s_j-s)|^2\right)\\
&=C_K\left(\int_{\mathcal C}|s_j-s|^2+\int_K|\nabla^{\mathcal E,k}s_j|^2\right)\\
&\longrightarrow 0.
\end{align*}
\end{proof}

\begin{proposition}\label{prop:parallel-section}
Suppose that $n$ is odd and $\mu\geq|J|_g$. Let $\pi:\widetilde M\to M$ be a connected regular covering, and let $p_j:M_j\to M$ be connected finite regular coverings admitting covering maps $q_j:\widetilde M\to M_j$ satisfying $p_j\circ q_j=\pi$. Suppose that there are $\varepsilon_j\longrightarrow0$ and smooth maps $N_j:\partial M_j\to S^{n-1}$ of positive degree such that
\[
\sqrt{(\tr_\Sigma k)^2+\|dN_j\|_{\tr}^2}\leq H_\Sigma+\varepsilon_j.
\]
Then the followings hold:
\begin{itemize}
	\item[(i)] The covering maps $q_j$ can be chosen so that a subsequence of $N_j\circ q_j$ converges locally uniformly to a smooth map $N:\partial\widetilde M\to S^{n-1}$ satisfying
\[
H_\Sigma=\sqrt{(\tr_\Sigma k)^2+\|dN\|_{\tr}^2}.
\]
	\item[(ii)] There is a smooth nowhere-zero section $s$ of $\mathcal E$ on $\widetilde M$ such that
\[
\nabla^{\mathcal E,k}s=0\quad\text{on }\widetilde M,\qquad \chi s=s\quad\text{on }\partial\widetilde M,
\]
where $\chi$ is defined by $N$.
\end{itemize}
\end{proposition}

\begin{proof}
Let $\mathcal E_j=p_j^*\mathcal E$ on $M_j$, and define $\chi_j=-\gamma(\nu)\otimes\gamma_0^*(N_j)$ on $\partial M_j$. By \cite[Proposition~2.4]{BBCH} applied to $M_j$ with boundary map $N_j$, we can find a non-trivial section $s_j$ of $\mathcal E_j$ such that $\left(\mathcal D-\frac{\ii}{2}\tr_g k\right)s_j=0$ on $M_j$, $\chi_js_j=s_j$ on $\partial M_j$, and $\int_{M_j}|s_j|^2=1$. By elliptic regularity, $s_j$ is smooth. The Sobolev trace inequality on $M$, pulled back to each cover, has a constant independent of $j$. Lemma~\ref{lem:sharp-boundary} and the dominant energy condition give
\[
\begin{aligned}
\int_{M_j}|\nabla^{\mathcal E,k}s_j|^2&\leq\frac{\varepsilon_j}{2}\int_{\partial M_j}|s_j|^2\\
&\leq C\varepsilon_j\left(\int_{M_j}|\nabla^{\mathcal E}s_j|^2+1\right)\\
&\leq C\varepsilon_j\left(\int_{M_j}|\nabla^{\mathcal E,k}s_j|^2+1\right).
\end{aligned}
\]
For sufficiently large $j$, this gives $\int_{M_j}|\nabla^{\mathcal E,k}s_j|^2\leq C\varepsilon_j$.

Fix a compact fundamental cell $\mathcal C$ for $\pi$, and choose increasing connected compact domains $K_\ell\supset\mathcal C$ exhausting $\widetilde M$. In the following sums, $q$ ranges over the covering maps $\widetilde M\to M_j$ satisfying $p_j\circ q=\pi$. Since $p_j$ is regular, there are $\deg p_j$ such maps, and their values at each $x\in\widetilde M$ are exactly the distinct points of $p_j^{-1}(\pi(x))$. Thus
\[
\sum_q\int_{\mathcal C}|q^*s_j|^2=\int_{M_j}|s_j|^2=1.
\]
For each $\ell$, the number of points in $K_\ell$ projecting to any fixed point of $M$ is bounded. Hence
\[
\begin{aligned}
\sum_q\int_{K_\ell}|\nabla^{\mathcal E,k}(q^*s_j)|^2&\leq C_\ell\int_{M_j}|\nabla^{\mathcal E,k}s_j|^2\\
&\leq C_\ell\varepsilon_j,
\end{aligned}
\]
where $C_\ell$ is independent of $j$. After passing to a subsequence and relabeling, we may assume $C_j\varepsilon_j\leq j^{-1}$. The two estimates above give a covering map $q_j$ with $\int_{\mathcal C}|q_j^*s_j|^2>0$ and
\[
\int_{K_j}|\nabla^{\mathcal E,k}(q_j^*s_j)|^2\leq\frac1j\int_{\mathcal C}|q_j^*s_j|^2.
\]

On $\widetilde M$, denote the normalized sections $q_j^*s_j/\bigl(\int_{\mathcal C}|q_j^*s_j|^2\bigr)^{1/2}$ again by $s_j$, and the boundary operators $q_j^*\chi_j=-\gamma(\nu)\otimes\gamma_0^*(N_j\circ q_j)$ again by $\chi_j$. Then $\chi_js_j=s_j$ on $\partial\widetilde M$, and
\[
\int_{\mathcal C}|s_j|^2=1,\qquad \int_{K_j}|\nabla^{\mathcal E,k}s_j|^2\leq\frac1j\longrightarrow0.
\]
Lemma~\ref{lem:compactness} then gives a subsequence converging strongly in $H^1_{\mathrm{loc}}$ up to the boundary to a smooth nowhere-zero section $s$ on $\widetilde M$ such that $\nabla^{\mathcal E,k}s=0$. The boundary traces converge strongly in $L^2_{\mathrm{loc}}(\partial\widetilde M,\mathcal E)$.

The assumption on $N_j$ gives a uniform Lipschitz bound for $N_j\circ q_j$ on $\partial\widetilde M$. After passing to a subsequence, these maps converge locally uniformly to a Lipschitz map $N:\partial\widetilde M\to S^{n-1}$. Let $\chi=-\gamma(\nu)\otimes\gamma_0^*(N)$ on $\partial\widetilde M$. Since $\chi_js_j=s_j$, for every compact $K\subset\partial\widetilde M$ we have
\[
\|\chi s-s\|_{L^2(K)}=\lim_{j\to\infty}\|\chi_js-s\|_{L^2(K)}\leq2\lim_{j\to\infty}\|s_j-s\|_{L^2(K)}=0.
\]
Thus $\chi s=s$, which implies $\gamma(\nu)s=s\gamma_0(N)$ on $\partial\widetilde M$. For an orthonormal basis $e_1^0,\ldots,e_n^0$ of $\R^n$, the Clifford relations give
\[
\begin{aligned}
\operatorname{Re}\langle\gamma(\nu)s,s\gamma_0(e_l^0)\rangle&=-\frac12\tr\left[s\bigl(\gamma_0(N)\gamma_0(e_l^0)+\gamma_0(e_l^0)\gamma_0(N)\bigr)s^*\right]\\
&=\langle N,e_l^0\rangle|s|^2.
\end{aligned}
\]
Since $s$ is smooth and nowhere zero on $\widetilde M$, each component $\langle N,e_l^0\rangle$ is smooth on $\partial\widetilde M$.

Fix a relatively compact open set $U\subset\partial\widetilde M$, and let $\varphi$ be a smooth compactly supported $\R^n$-valued one-form on $U$ with $\|\varphi\|_{\operatorname{op}}\leq1$. Integration by parts and the local uniform convergence of $N_j\circ q_j$ give
\[
\begin{aligned}
\int_U\langle dN,\varphi\rangle&=\lim_{j\to\infty}\int_U\langle d(N_j\circ q_j),\varphi\rangle\\
&\leq\lim_{j\to\infty}\int_U\sqrt{(H_\Sigma+\varepsilon_j)^2-(\tr_\Sigma k)^2}\\
&=\int_U\sqrt{H_\Sigma^2-(\tr_\Sigma k)^2}.
\end{aligned}
\]
This gives
\[
\int_U\|dN\|_{\tr}=\sup_\varphi\int_U\langle dN,\varphi\rangle\leq\int_U\sqrt{H_\Sigma^2-(\tr_\Sigma k)^2}.
\]
Since $U$ is arbitrary and $N$ is smooth, $H_\Sigma\geq\sqrt{(\tr_\Sigma k)^2+\|dN\|_{\tr}^2}$ on $\partial\widetilde M$.

Let $e_1,\ldots,e_{n-1}$ be a local orthonormal frame tangent to $\partial\widetilde M$. Since $\nabla^{\mathcal E,k}s=0$ on $\widetilde M$ and $\gamma(\nu)s=s\gamma_0(N)$ on $\partial\widetilde M$, we have
\[
\begin{aligned}
\sum_a\gamma(e_a)s\gamma_0(dN(e_a))&=\sum_a\gamma(e_a)\nabla_{e_a}^{\mathcal E,k}(s\gamma_0(N))\\
&=\sum_a\gamma(e_a)\nabla_{e_a}^{\mathcal E,k}(\gamma(\nu)s)\\
&=\sum_a\gamma(e_a)\left[\gamma(D_{e_a}\nu)+\frac{\ii}{2}\bigl(\gamma(k(e_a))\gamma(\nu)-\gamma(\nu)\gamma(k(e_a))\bigr)\right]s\\
&=\sum_{a,b}h(e_a,e_b)\gamma(e_a)\gamma(e_b)s+\ii\sum_{a,b}k(e_a,e_b)\gamma(e_a)\gamma(e_b)\gamma(\nu)s\\
&=-H_\Sigma s-\ii(\tr_\Sigma k)\gamma(\nu)s.
\end{aligned}
\]
By \eqref{eq:boundary-operator-bound},
\[
\begin{aligned}
H_\Sigma^2|s|^2&=\left|\sum_a\gamma(e_a)s\gamma_0(dN(e_a))+\ii(\tr_\Sigma k)\gamma(\nu)s\right|^2\\
&\leq\bigl((\tr_\Sigma k)^2+\|dN\|_{\tr}^2\bigr)|s|^2\\
&\leq H_\Sigma^2|s|^2.
\end{aligned}
\]
Since $s$ is nowhere zero, this proves the asserted equality.
\end{proof}

\begin{corollary}\label{cor:parallel-frame}
Suppose that $\mu\geq|J|_g$. Suppose that there are $\varepsilon_j\longrightarrow0$ and smooth maps $N_j:\Sigma\to S^{n-1}$ of positive degree such that
\[
\sqrt{(\tr_\Sigma k)^2+\|dN_j\|_{\tr}^2}\leq H_\Sigma+\varepsilon_j.
\]
Then the followings hold:
\begin{itemize}
	\item[(i)] A subsequence of $N_j$ converges uniformly to a smooth map $N$ of positive degree satisfying
\[
H_\Sigma=\sqrt{(\tr_\Sigma k)^2+\|dN\|_{\tr}^2}.
\]
	\item[(ii)] There is a smooth section $s$ of $\mathcal E$ on $M$ such that
\[
\nabla^{\mathcal E,k}s=0\quad\text{on }M,\qquad \gamma(\nu)\circ s=s\circ\gamma_0(N)\quad\text{on }\Sigma.
\]
	\item[(iii)] The section $s$, considered as a homomorphism field $s:S_0\to\mathcal S$, is invertible at each point in $M$.
\end{itemize} 
\end{corollary}

\begin{proof}
If $n$ is odd, apply Proposition~\ref{prop:parallel-section} with $\widetilde M=M_j=M$ and all covering maps equal to the identity. We obtain uniform convergence of a subsequence of $N_j$ to a smooth map $N$, the asserted equality, and a smooth nowhere-zero section $s$ with $\nabla^{\mathcal E,k}s=0$ on $M$ and $\gamma(\nu)s=s\gamma_0(N)$ on $\Sigma$. For large $j$, the maps $N_j$ and $N$ are homotopic, so $\deg N>0$.

Suppose now that $n$ is even. The hypothesis bounds the Lipschitz constants of $N_j$, so a subsequence converges uniformly to a Lipschitz map $N$. Again $\deg N>0$. The integral argument in the proof of Proposition~\ref{prop:parallel-section} gives
\begin{equation}\label{eq:limit-boundary-comparison}
\|dN\|_{\tr}^2\leq H_\Sigma^2-(\tr_\Sigma k)^2
\end{equation}
almost everywhere on on $\Sigma$. Consider $M\times S^1$ with the metric $g+dt^2$. For orthonormal frames $e_1,\ldots,e_n$ and $e_1^0,\ldots,e_n^0$, extend the Clifford multiplications by
\[
\gamma(\partial_t)=\ii^{n/2+1}\gamma(e_1)\cdots\gamma(e_n),\qquad \gamma_0(e_{n+1}^0)=\ii^{n/2+1}\gamma_0(e_1^0)\cdots\gamma_0(e_n^0).
\]
We consider the coverings
\[
M\times\R\longrightarrow M\times(\R/2\pi j\Z)\xrightarrow{p_j}M\times(\R/2\pi\Z),
\]
where both maps are the natural projections and $p_j$ has degree $j$. Denote their composite by $\pi$. Let $\Phi:S^{n-1}\times S^1\to S^n$ be the degree-one map in \cite[Lemma~3.8]{BBCH}, and define on $\Sigma\times(\R/2\pi j\Z)$
\[
\widetilde N_j(x,t)=\Phi(N_j(x),t/j).
\]
Then $\deg\widetilde N_j=\deg N_j>0$, and \cite[Lemma~3.9]{BBCH} gives
\[
\begin{aligned}
\sqrt{(\tr_\Sigma k)^2+\|d\widetilde N_j\|_{\tr}^2}&\leq\sqrt{(\tr_\Sigma k)^2+(\|dN_j\|_{\tr}+2/j)^2}\\
&\leq\sqrt{(\tr_\Sigma k)^2+\|dN_j\|_{\tr}^2}+\frac2j\\
&\leq H_\Sigma+\varepsilon_j+\frac2j.
\end{aligned}
\]
Proposition~\ref{prop:parallel-section} therefore applies with $\varepsilon_j+2/j\to0$. It gives covering maps $q_j:M\times\R\to M\times(\R/2\pi j\Z)$, local uniform convergence of $\widetilde N_j\circ q_j$ to a smooth map $\widetilde N:\Sigma\times\R\to S^n$, and a smooth nowhere-zero parallel section $s$ on $M\times\R$ satisfying
\[
H_\Sigma^2-(\tr_\Sigma k)^2=\|d\widetilde N\|_{\tr}^2,\qquad \gamma(\nu)s=s\gamma_0(\widetilde N)\quad\text{on }\Sigma\times\R.
\]
Write $q_j(x,t)=(x,\tau_j(x,t))$, where $\tau_j$ takes values in $\R/2\pi j\Z$. Since $p_j\circ q_j=\pi$,
\[
\left(e^{\ii(\tau_j(x,t)-t)/j}\right)^j=e^{\ii\tau_j(x,t)}e^{-\ii t}=1.
\]
The continuous map $e^{\ii(\tau_j-t)/j}$ takes values in the finite set of $j$-th roots of unity. Since $M\times\R$ is connected, this map is constant. Thus there is $\theta_j\in S^1$ such that $\tau_j(x,t)/j=\theta_j+t/j$ in $S^1$ for every $(x,t)\in M\times\R$. Pass to a subsequence we get $\theta_j\to\theta_*$. Since $N_j\to N$ uniformly and $t/j\to0$ uniformly on bounded intervals,
\[
\begin{aligned}
\widetilde N(x,t)&=\lim_{j\to\infty}\Phi\left(N_j(x),\frac{\tau_j(x,t)}j\right)\\
&=\lim_{j\to\infty}\Phi\left(N_j(x),\theta_j+\frac tj\right)\\
&=\Phi(N(x),\theta_*).
\end{aligned}
\]
If $\Phi(\cdot,\theta_*)$ has Lipschitz constant $a<1$, then \eqref{eq:limit-boundary-comparison} gives
\[
\|dN\|_{\tr}^2\leq H_\Sigma^2-(\tr_\Sigma k)^2=\|d\widetilde N\|_{\tr}^2\leq a^2\|dN\|_{\tr}^2
\]
almost everywhere. Hence $N$ is locally constant, contrary to $\deg N>0$. By the definition of $\Phi$ in \cite[Lemma~3.8]{BBCH}, we have $\Phi(\cdot,\theta_*) = (\cdot, 0)$ which is the standard equatorial inclusion. Thus $\widetilde N=(N,0)$, $N$ is smooth, and
\[
H_\Sigma^2-(\tr_\Sigma k)^2=\|d\widetilde N\|_{\tr}^2=\|dN\|_{\tr}^2.
\]
Restricting to $M\times\{0\}$ gives a smooth nowhere-zero section $s$ on $M$ with $\nabla^{\mathcal E,k}s=0$ and $\gamma(\nu)s=s\gamma_0(N)$ on $\Sigma$.

We next prove invertibility of $s$ as in \cite[Lemmas~4.3--4.5]{BBCH}. Define the linear subspace
\[
L=\{\xi\in S_0:s_x\xi=0\text{ for every }x\in M\}.
\]
For each fixed $\xi\in S_0$,
\[
\nabla_X^{\mathcal S}(s\xi)+\frac{\ii}{2}\gamma(k(X))s\xi=(\nabla_X^{\mathcal E,k}s)\xi=0.
\]
If $s_x\xi=0$ at one point $x\in M$, ODE uniqueness and connectedness of $M$ give $s_x\xi=0$ for every $x\in M$. Consequently, $\ker s_x=L$ for every $x\in M$.

For $p\in\Sigma$ and $\xi\in L$, the boundary condition yields
\[
s_p\gamma_0(N(p))\xi=\gamma(\nu(p))s_p\xi=0\quad\implies\quad \gamma_0(N(p))\xi\in\ker s_p=L.
\]
Since $\deg N>0$, we have $N(\Sigma)=S^{n-1}$. It follows that $\gamma_0(v)L\subset L$ for every $v\in\R^n$. These endomorphisms generate $\operatorname{End}(S_0)$, hence $L=\{0\}$ or $L=S_0$. The latter would give $s=0$, so $\ker s_x=L=\{0\}$ at every $x\in M$. Since $S_0$ and $\mathcal S_x$ have the same complex dimension, $s_x$ is invertible.
\end{proof}

\section{Parallel spinors and rigidity}\label{sec:comparison-metric}

In this section, let $(M^n,g,k)$ be a connected spin initial data set with boundary $\Sigma$, and we use the notations of the previous section. Let $N:\Sigma\to S^{n-1}$ be smooth, and let $s$ be a nonzero smooth section of $\mathcal E$ satisfying
\[
\nabla^{\mathcal E,k}s=0\quad\text{on }M,\qquad \chi s=s\quad\text{on }\Sigma.
\]
Define the vector field $V$ and its $g$-dual $V^\flat$ by
\begin{equation}\label{eq:current-definition}
V^\flat(Y)=g(V,Y)=\operatorname{Re}\langle(\ii\gamma(Y)\otimes\Id)s,s\rangle.
\end{equation}

\begin{lemma}\label{lem:current}
For every vector field $X$, we have
\[
D_XV=-|s|^2k(X),\qquad X(|s|^2)=-k(X,V).
\]
Moreover, $dV^\flat=0$, $|s|^4-|V|_g^2$ is constant, and $V^\flat|_{T\Sigma}=0$.
\end{lemma}

\begin{proof}
Fix an orthonormal basis $\{\xi_\alpha\}$ of $S_0$. Since $\nabla_X^{\mathcal S}(s\xi_\alpha)+\frac{\ii}{2}\gamma(k(X))s\xi_\alpha=0$, applying \cite[Lemma~4.16]{BBCH} to $u=v=s\xi_\alpha$ and summing over $\alpha$ gives the first two identities.
Since $k$ is symmetric, $dV^\flat=0$. Furthermore,
\[
\begin{aligned}
X(|s|^4-|V|_g^2)&=2|s|^2X(|s|^2)-2g(D_XV,V)\\
&=-2|s|^2k(X,V)+2|s|^2k(X,V)\\
&=0.
\end{aligned}
\]
The assertion $V^\flat|_{T\Sigma}=0$ follows from \eqref{eq:tangential-current}.
\end{proof}

\begin{proposition}\label{prop:comparison-metric}
Suppose that $|s|^4-|V|_g^2=1$, and define the metric
\begin{equation}\label{eq:comparison-metric}
\bar g=g+V^\flat\otimes V^\flat.
\end{equation}
Then $\bar g$ is Ricci flat and $\bar g|_\Sigma=g|_\Sigma$ on $\Sigma$. The second fundamental form $\bar h$ and mean curvature $\bar H_\Sigma$ of $\Sigma$ with respect to $\bar g$ satisfy
\[
\bar h=\frac{h}{|s|^2},\qquad \bar H_\Sigma=\frac{H_\Sigma}{|s|^2},\qquad \bar H_\Sigma^2=H_\Sigma^2-(\tr_\Sigma k)^2.
\]
Moreover, if $X\in T\Sigma$ satisfies $dN(X)=0$, then $\bar h(X,\cdot)=0$.
\end{proposition}

\begin{proof}
Write $\bar D$ for the Levi--Civita connection of $\bar g$. Lemma~\ref{lem:current} gives
\[
(D_X\bar g)(Y,Z)=-|s|^2\bigl(k(X,Y)V^\flat(Z)+k(X,Z)V^\flat(Y)\bigr).
\]
This implies
\[
\begin{aligned}
2\bar g(\bar D_XY-D_XY,Z)&=(D_X\bar g)(Y,Z)+(D_Y\bar g)(X,Z)-(D_Z\bar g)(X,Y)\\
&=-2|s|^2k(X,Y)V^\flat(Z).
\end{aligned}
\]
Since $\bar g(V,Z)=(1+|V|_g^2)V^\flat(Z)=|s|^4V^\flat(Z)$, it follows that
\begin{equation}\label{eq:comparison-connection}
\bar D_XY=D_XY-\frac{k(X,Y)}{|s|^2}V.
\end{equation}

On a simply connected coordinate neighborhood $U$, we can write $df=-V^\flat$. Lemma~\ref{lem:current} gives $\Hess_g f=|s|^2k$, and
\[
1-|df|_{\bar g}^2=1-\frac{|V|_g^2}{1+|V|_g^2}=|s|^{-4}.
\]
The graph $x\mapsto(x,f(x))$ in $(U\times\R,\bar g-dt^2)$ has induced metric $g$. By \eqref{eq:comparison-connection}, its second fundamental form is
\[
\frac{\Hess_{\bar g}f}{\sqrt{1-|df|_{\bar g}^2}}=|s|^2\left(\Hess_g f-\frac{|V|_g^2}{|s|^2}k\right)=(|s|^4-|V|_g^2)k=k.
\]
We now consider the product manifold $Z=U\times\R$ with metric $\bar g-dt^2$. We write $T=\partial_t$ and $\widehat X=X+df(X)T$ for $X\in TU$. Let $\xi\in S_0$ be such that $\psi=s\xi$ is nonzero at one point of $U$. Since $\nabla^{\mathcal S,k}\psi=0$, ODE uniqueness implies that $\psi$ is nowhere zero on $U$. Along the graph $x\mapsto(x,f(x))$, the hypersurface spinor identification is given by
\[
\gamma(X)\psi=\ii\nu_Z\cdot\widehat X\cdot\psi,\qquad \nu_Z=\frac{T+\nabla^{\bar g}f}{\sqrt{1-|df|_{\bar g}^2}},
\]
where $\cdot$ denotes ambient Clifford multiplication. For odd $n$, $\mathcal S$ is identified with the appropriate ambient half-spin bundle. The spin connection formula \cite[(2.5)]{BGM} gives
\[
\begin{aligned}
\nabla^{\mathcal S_Z}_{\widehat X}\psi&=\nabla_X^{\mathcal S}\psi-\frac12\nu_Z\cdot\widehat{k(X)}\cdot\psi\\
&=\nabla_X^{\mathcal S}\psi+\frac{\ii}{2}\gamma(k(X))\psi\\
&=0.
\end{aligned}
\]
Consequently, along the graph $x\mapsto(x,f(x))$, we have
\[
\begin{aligned}
R^{\mathcal S_Z}(\widehat X,\widehat Y)\psi&=\left(\nabla^{\mathcal S_Z}_{\widehat X}\nabla^{\mathcal S_Z}_{\widehat Y}-\nabla^{\mathcal S_Z}_{\widehat Y}\nabla^{\mathcal S_Z}_{\widehat X}-\nabla^{\mathcal S_Z}_{\widehat{[X,Y]}}\right)\psi=0.
\end{aligned}
\]
Since $Z$ is a product, $R^{\mathcal S_Z}(T,\cdot)=0$. Thus, regarding $X,Y$ also as horizontal vectors, we have
\[
R^{\mathcal S_Z}(X,Y)\psi=R^{\mathcal S_Z}(\widehat X,\widehat Y)\psi=0.
\]
Let $E_1,\ldots,E_n$ be a $\bar g$-orthonormal frame, regarded as horizontal vectors in $TZ$. By \cite[(2.6)]{BGM},
\[
\begin{aligned}
\frac12\Ric_{\bar g}(X)\cdot\psi&=\sum_{a=1}^nE_a\cdot R^{\mathcal S_Z}(E_a,X)\psi-T\cdot R^{\mathcal S_Z}(T,X)\psi=0.
\end{aligned}
\]
It follows that
\[
0=\Ric_{\bar g}(X)\cdot\Ric_{\bar g}(X)\cdot\psi=-|\Ric_{\bar g}(X)|_{\bar g}^2\psi,
\]
and hence $\Ric_{\bar g}=0$.

On $\Sigma$, Lemma~\ref{lem:current} gives $V=g(V,\nu)\nu$, so $\bar g|_\Sigma=g|_\Sigma$ and $\bar\nu=\nu/|s|^2$. For $X,Y\in T\Sigma$, \eqref{eq:comparison-connection} yields
\[
\begin{aligned}
\bar h(X,Y)&=\bar g\left(\bar D_X\frac{\nu}{|s|^2},Y\right)\\
&=\frac1{|s|^2}\bar g\left(D_X\nu-\frac{k(X,\nu)}{|s|^2}V,Y\right)+X(|s|^{-2})\bar g(\nu,Y)\\
&=\frac{h(X,Y)}{|s|^2}.
\end{aligned}
\]
Also,
\[
-|s|^2k(X,Y)=g(D_XV,Y)=g(V,\nu)h(X,Y).
\]
This gives $\bar H_\Sigma=H_\Sigma/|s|^2$ and $\tr_\Sigma k=-g(V,\nu)H_\Sigma/|s|^2$. Consequently,
\[
\begin{aligned}
H_\Sigma^2-(\tr_\Sigma k)^2&=H_\Sigma^2\left(1-\frac{g(V,\nu)^2}{|s|^4}\right)
=\frac{H_\Sigma^2}{|s|^4}=\bar H_\Sigma^2.
\end{aligned}
\]

Let $X\in T\Sigma$ satisfy $dN(X)=0$, and write $h(X)=D_X\nu$. For a local orthonormal frame $e_1,\ldots,e_{n-1}$ tangent to $\Sigma$, the equations $\nabla^{\mathcal E,k}s=0$ and $\gamma(\nu)s=s\gamma_0(N)$ on $\Sigma$ give
\[
\begin{aligned}
0=s\gamma_0(dN(X))&=\nabla_X^{\mathcal E,k}(\gamma(\nu)s)\\
&=\gamma(h(X))s+\frac{\ii}{2}\bigl(\gamma(k(X))\gamma(\nu)-\gamma(\nu)\gamma(k(X))\bigr)s\\
&=\gamma(h(X))s+\ii\sum_a k(X,e_a)\gamma(e_a)\gamma(\nu)s\\
&=\gamma(h(X))\left(\Id-\frac{\ii g(V,\nu)}{|s|^2}\gamma(\nu)\right)s.
\end{aligned}
\]
Since $(\ii\gamma(\nu))^2=\Id$ and $|g(V,\nu)|<|s|^2$, the endomorphism $\Id-\frac{\ii g(V,\nu)}{|s|^2}\gamma(\nu)$ is invertible. If $h(X)\ne0$, then $\gamma(h(X))$ is also invertible, contradicting $s\ne0$. Thus $h(X,\cdot)=0$, and consequently $\bar h(X,\cdot)=0$.
\end{proof}

\begin{lemma}\label{lem:metric-change}
Let $(M,g)$ be a Riemannian manifold with boundary $\Sigma$, and let $f\in C^\infty(M)$ be constant on $\Sigma$. Suppose that $\bar g=g+df\otimes df$ is flat and $k=(1+|df|_g^2)^{-1/2}\Hess_g f$. Then $\mu=J=0$, $\bar g_\Sigma = g_\Sigma$, and
\[
H_\Sigma=\sqrt{1+\nu(f)^2}\,\bar H_\Sigma,\qquad \tr_\Sigma k=\nu(f)\bar H_\Sigma,
\]
where $\nu$ is the outward unit normal for $g$ and $\bar H_\Sigma$ is the boundary mean curvature for $\bar g$.
\end{lemma}

\begin{proof}
Write $D$ and $\bar D$ for the Levi--Civita connections of $g$ and $\bar g$. we have
\[
\begin{aligned}
2\bar g(\bar D_XY-D_XY,Z)&=(D_X\bar g)(Y,Z)+(D_Y\bar g)(X,Z)-(D_Z\bar g)(X,Y)\\
&=2\Hess_g f(X,Y)df(Z),
\end{aligned}
\]
and hence $\bar D_XY=D_XY+(1+|df|_g^2)^{-1}\Hess_g f(X,Y)\nabla^g f$. Since $f$ is constant on $\Sigma$, we have $\bar g_\Sigma = g_\Sigma$ and $\bar\nu=(1+\nu(f)^2)^{-1/2}\nu$. For $X,Y\in T\Sigma$, we obtain
\[
\bar h(X,Y)=\bar g(\bar D_X\bar\nu,Y)=\frac{g(D_X\nu,Y)}{\sqrt{1+\nu(f)^2}}=\frac{h(X,Y)}{\sqrt{1+\nu(f)^2}}
\]
and
\[
\begin{aligned}
k(X,Y)&=\frac{\Hess_g f(X,Y)}{\sqrt{1+\nu(f)^2}}=\frac{-df(D_XY)}{\sqrt{1+\nu(f)^2}}
=\frac{\nu(f)h(X,Y)}{\sqrt{1+\nu(f)^2}}=\nu(f)\bar h(X,Y).
\end{aligned}
\]
Taking traces proves the asserted identities. In local coordinates for $\bar g$, the remark following Theorem~\ref{thm:boundary-rigidity} identifies $(g,k)$ with the induced metric and second fundamental form of a spacelike graph. The Gauss and Codazzi equations give
\[
R_g=|k|_g^2-(\tr_g k)^2,\qquad \operatorname{div}_g k=d\tr_g k,
\]
so $\mu=J=0$.
\end{proof}

\subsection{Proof of Theorem~\ref{thm:boundary-rigidity} and Corollary~\ref{cor:normalized-comparison}}
\begin{proof}[Proof of Theorem~\ref{thm:boundary-rigidity}]
Apply Corollary~\ref{cor:parallel-frame} with $N_j=N$ and $\varepsilon_j=0$. We obtain the equality $H_\Sigma=\sqrt{(\tr_\Sigma k)^2+\|dN\|_{\tr}^2}$ and an everywhere invertible section $s$ with $\nabla^{\mathcal E,k}s=0$ and $\chi s=s$ on $\Sigma$. Define $V$ by \eqref{eq:current-definition}. If $V\equiv0$, Lemma~\ref{lem:current} allows us to rescale $s$ so that $|s|^4-|V|_g^2=1$. Otherwise, at a point where $V\neq0$, let $X=V/|V|_g$ and let $Y\perp X$ be unit. Since $\ii\gamma(X)$ anticommutes with the invertible $\gamma(Y)$, we have $\Id-\ii\gamma(X)\neq0$. Since $s$ is onto,
\[
|s|^2-|V|_g=\frac12\left|((\Id-\ii\gamma(X))\otimes\Id)s\right|^2>0.
\]
Thus $|s|^4-|V|_g^2$ is a positive constant by Lemma~\ref{lem:current}, and we again rescale $s$ so that it equals $1$. Define $\bar g$ by \eqref{eq:comparison-metric}.

Write $\bar D$ and $\bar R$ for the Levi--Civita connection and Riemann curvature tensor of $(M,\bar g)$. For $Y\in TM$, define
\[
\bar\gamma(Y)=s^{-1}\left(|s|^2\gamma(Y)+\frac{\ii}{2}[\gamma(V),\gamma(Y)]\right)s\in\operatorname{End}(S_0).
\]
The Clifford relations give
\[
\begin{aligned}
\bar\gamma(Y)^2&=-\bigl(|s|^4|Y|_g^2-|V|_g^2|Y|_g^2+g(V,Y)^2\bigr)\Id
=-\bar g(Y,Y)\Id,
\end{aligned}
\]
so $\bar\gamma$ is injective. By Lemma~\ref{lem:current}, \eqref{eq:comparison-connection}, and $\nabla^{\mathcal E,k}s=0$, we obtain
\[
\begin{aligned}
X(\bar\gamma(Y))&=\bar\gamma(D_XY)+s^{-1}\left(-k(X,V)\gamma(Y)-\frac14[\gamma(k(X)),[\gamma(V),\gamma(Y)]]\right)s\\
&=\bar\gamma(D_XY)-k(X,Y)s^{-1}\gamma(V)s\\
&=\bar\gamma\left(D_XY-\frac{k(X,Y)}{|s|^2}V\right)\\
&=\bar\gamma(\bar D_XY).
\end{aligned}
\]
Consequently,
\[
\begin{aligned}
\bar\gamma(\bar R(X,Y)Z)&=X\bigl(Y(\bar\gamma(Z))\bigr)-Y\bigl(X(\bar\gamma(Z))\bigr)-[X,Y](\bar\gamma(Z))=0.
\end{aligned}
\]
Injectivity of $\bar\gamma$ gives $\bar R=0$ on $M$.

Along $\Sigma$, we have $V=g(V,\nu)\nu$ and $\bar\nu=\nu/|s|^2$, so the chiral boundary condition gives
\[
\bar\gamma(\bar\nu)=s^{-1}\gamma(\nu)s=\gamma_0(N).
\]
For $X\in T\Sigma$, it follows that
\[
\begin{aligned}
|dN(X)|^2\Id&=-\gamma_0(dN(X))^2=-\bigl[X(\bar\gamma(\bar\nu))\bigr]^2
=-\bar\gamma(\bar D_X\bar\nu)^2=|\bar D_X\bar\nu|_{\bar g}^2\Id.
\end{aligned}
\]
Since $\bar g|_\Sigma=g|_\Sigma$, the singular values of $dN$ are the absolute values of the eigenvalues $\kappa_1,\ldots,\kappa_{n-1}$ of $\bar h$. Since $H_\Sigma\geq0$, Proposition~\ref{prop:comparison-metric} and Corollary~\ref{cor:parallel-frame} imply
\[
\sum_a\kappa_a=\bar H_\Sigma=\sqrt{H_\Sigma^2-(\tr_\Sigma k)^2}=\|dN\|_{\tr}=\sum_a|\kappa_a|.
\]
Hence $\bar h\geq0$. Since $\deg N>0$, there is a point where $dN$ has full rank, and $\bar h$ is positive definite there. By \cite[Theorem~A.1]{BBCH}, $(M,\bar g)$ is isometric to a compact convex domain $K\subset\R^n$ with smooth boundary.

The domain $K$ is contractible, so Lemma~\ref{lem:current} gives a smooth function $f$ with $df=-V^\flat$. Since $\partial K$ is connected and $V^\flat|_{T\Sigma}=0$, the function $f$ is constant on the boundary. We have $\bar g=g+df\otimes df$ and
\[
(1+|df|_g^2)^{-1/2}\Hess_g f=\frac{-DV^\flat}{\sqrt{1+|V|_g^2}}=\frac{|s|^2k}{|s|^2}=k.
\]

Identify $(M,\bar g)$ with $K$, and let $N_0$ be the outward Gauss map of $\partial K$. For a constant vector $v\in\R^n$, the identity $X(\bar\gamma(v))=\bar\gamma(\bar D_Xv)=0$ shows that $\bar\gamma$ is a fixed linear map from $\R^n$ to $\operatorname{End}(S_0)$. Along $\partial K$, we have $\bar\gamma(N_0)=\gamma_0(N)$. Since $N$ is onto, $\gamma_0(\R^n)\subseteq\bar\gamma(\R^n)$. Both spaces have real dimension $n$, so there is a linear map $O:\R^n\to\R^n$ satisfying $\bar\gamma(Ov)=\gamma_0(v)$. Moreover,
\[
|Ov|^2\Id=-\bar\gamma(Ov)^2=-\gamma_0(v)^2=|v|^2\Id.
\]
Thus $O\in O(n)$, and $\bar\gamma(N_0)=\gamma_0(N)=\bar\gamma(ON)$ implies $N_0=ON$.

Conversely, suppose that $f$ is constant on $\Sigma$, $k=(1+|df|_g^2)^{-1/2}\Hess_g f$, and $(M,\bar g)$ is isometric to a smooth compact convex domain $K$. Let $N_0$ be its outward Gauss map. By convexity, $\bar H_\Sigma=\|dN_0\|_{\tr}$, and Lemma~\ref{lem:metric-change} gives $\mu=J=0$, $H_\Sigma\geq0$, and
\[
\begin{aligned}
H_\Sigma^2-(\tr_\Sigma k)^2&=(1+\nu(f)^2)\bar H_\Sigma^2-\nu(f)^2\bar H_\Sigma^2
=\bar H_\Sigma^2=\|dN_0\|_{\tr}^2.
\end{aligned}
\]
The domain $K$ is spin, and $\deg N_0=1$, proving the converse.	
\end{proof}

\begin{proof}[Proof of Corollary~\ref{cor:normalized-comparison}]
The tensor inequality and $H_\Sigma\geq0$ give $H_\Sigma\geq|\tr_\Sigma k|$. Let $N_0$ be the outward Gauss map of $\partial K_0$. For $x\in\Sigma$, let $e_1^0,\ldots,e_{n-1}^0$ be an orthonormal eigenbasis of $dN_0$ at $F(x)$ with eigenvalues $\kappa_a^0>0$. The assumption gives
\[
\begin{aligned}
\|d(N_0\circ F)\|_{\tr}&=\left\|\sum_a\kappa_a^0 e_a^0\otimes((dF)^*e_a^0)^\flat\right\|_{\tr}\\
&\leq\sum_a\kappa_a^0|(dF)^*e_a^0|\\
&\leq\frac{\sqrt{H_\Sigma^2-(\tr_\Sigma k)^2}}{H_0\circ F}\sum_a\kappa_a^0\\
&=\sqrt{H_\Sigma^2-(\tr_\Sigma k)^2}.
\end{aligned}
\]
Since $\deg(N_0\circ F)=\deg F>0$, Theorem~\ref{thm:boundary-rigidity} gives a smooth function $f$, with $f$ constant on $\Sigma$ and $k=(1+|df|_g^2)^{-1/2}\Hess_g f$, such that $(M,\bar g)$ is isometric to a compact convex Euclidean domain $K$, where $\bar g=g+df\otimes df$. Identifying $(M,\bar g)$ with $K$ and composing this identification with an orthogonal transformation, we have $N_K=N_0\circ F$.

We then have $T_x\partial K=T_{F(x)}\partial K_0$ as subspaces of $\R^n$. Lemma~\ref{lem:metric-change} gives $\bar H_\Sigma=\sqrt{H_\Sigma^2-(\tr_\Sigma k)^2}$. This yields
\[
\begin{aligned}
\bar H_\Sigma&=\tr(dN_K)=\sum_a\kappa_a^0\langle dF(e_a^0),e_a^0\rangle
\leq\sum_a\kappa_a^0|dF(e_a^0)|\leq\frac{\bar H_\Sigma}{H_0\circ F}\sum_a\kappa_a^0=\bar H_\Sigma.
\end{aligned}
\]
Since every $\kappa_a^0$ is positive, equality holds in each summand, so
\[
\langle dF(e_a^0),e_a^0\rangle=|dF(e_a^0)|=\frac{\bar H_\Sigma}{H_0\circ F},\qquad dF(e_a^0)=\frac{\bar H_\Sigma}{H_0\circ F}e_a^0.
\]
Let $x:\partial K\hookrightarrow\R^n$ denote the inclusion. Then $dF=(\bar H_\Sigma/(H_0\circ F))dx$. For each unit tangent vector $X$, let $Y$ be a unit tangent vector orthogonal to $X$. Since $d(dF)=0$,
\[
\begin{aligned}
0&=\langle d(dF)(X,Y),Y\rangle\\
&=\left\langle X\left(\frac{\bar H_\Sigma}{H_0\circ F}\right)Y-Y\left(\frac{\bar H_\Sigma}{H_0\circ F}\right)X,Y\right\rangle\\
&=X\left(\frac{\bar H_\Sigma}{H_0\circ F}\right).
\end{aligned}
\]
Thus $\bar H_\Sigma/(H_0\circ F)=a\geq0$, and connectedness gives $F(x)=ax+b$ for some $b\in\R^n$. Since $\deg F>0$, the map $F$ is onto $\partial K_0$ and $a>0$. Hence $K_0=aK+b$, so $F$ extends to an isometry from $(M,a^2\bar g)$ onto $K_0$.

Conversely, suppose that $f$ and $a$ satisfy the asserted conclusions and that $F$ extends to an  isometry from $(M,a^2\bar g)$ onto $K_0$. Since $g$ and $\bar g$ agree on $T\Sigma$, the isometry and Lemma~\ref{lem:metric-change} give
\[
\begin{aligned}
\bigl(H_\Sigma^2-(\tr_\Sigma k)^2\bigr)g_\Sigma&=\bar H_\Sigma^2g_\Sigma
=(H_0\circ F)^2a^2g_\Sigma=(H_0\circ F)^2F^*g_{\R^n}.
\end{aligned}
\]
The remaining hypotheses follow from the converse of Theorem~\ref{thm:boundary-rigidity}.	
\end{proof}

\section{The hyperspherical radius}\label{sec:radius}

\begin{proof}[Proof of Theorem~\ref{thm:radius}]
Suppose that
\[
H_\Sigma\geq\sqrt{(\tr_\Sigma k)^2+\frac{(n-1)^2}{\operatorname{Rad}(\Sigma)^2}}\qquad\text{on }\Sigma.
\]
By the definition of $\operatorname{Rad}(\Sigma)$, we can find smooth maps $N_j:\Sigma\to S^{n-1}$ such that
\[
\deg N_j>0,\qquad \operatorname{Lip}(N_j)\leq\frac1{\operatorname{Rad}(\Sigma)}+\frac1j,
\]
Since $\|dN_j\|_{\tr}\leq(n-1)(\operatorname{Rad}(\Sigma)^{-1}+j^{-1})$, the triangle inequality gives
\[
\begin{aligned}
\sqrt{(\tr_\Sigma k)^2+\|dN_j\|_{\tr}^2}&\leq\sqrt{(\tr_\Sigma k)^2+\left(\frac{n-1}{\operatorname{Rad}(\Sigma)}+\frac{n-1}{j}\right)^2}\\
&\leq\sqrt{(\tr_\Sigma k)^2+\frac{(n-1)^2}{\operatorname{Rad}(\Sigma)^2}}+\frac{n-1}{j}\\
&\leq H_\Sigma+\frac{n-1}{j}.
\end{aligned}
\]
Corollary~\ref{cor:parallel-frame}, with $\varepsilon_j=(n-1)/j$, gives a subsequence converging uniformly to a smooth map $N:\Sigma\to S^{n-1}$ of positive degree satisfying
\[
H_\Sigma=\sqrt{(\tr_\Sigma k)^2+\|dN\|_{\tr}^2}.
\]
For $x,y\in\Sigma$, uniform convergence gives
\[
\begin{aligned}
d_{S^{n-1}}(N(x),N(y))&=\lim_j d_{S^{n-1}}(N_j(x),N_j(y))\\
&\leq\lim_j\left(\frac1{\operatorname{Rad}(\Sigma)}+\frac1j\right)d_\Sigma(x,y)\\
&=\frac{d_\Sigma(x,y)}{\operatorname{Rad}(\Sigma)}.
\end{aligned}
\]
Thus $\operatorname{Lip}(N)\leq\operatorname{Rad}(\Sigma)^{-1}$.

By Theorem~\ref{thm:boundary-rigidity}, there is $f\in C^\infty(M)$, such that $f$ is constant on $\Sigma$, $k=(1+|df|_g^2)^{-1/2}\Hess_g f$ and $(M,\bar g)$ is isometric to a compact convex Euclidean domain $K$ with smooth boundary, where $\bar g=g+df\otimes df$. Its outward Gauss map $N_0$ agrees with $N$ up to an orthogonal transformation. Since $\bar g_\Sigma = g_\Sigma$, we have $\operatorname{Lip}(N_0)\leq\operatorname{Rad}(\Sigma)^{-1}$. Let $\kappa_1,\ldots,\kappa_{n-1}$ be the principal curvatures of $\partial K$. Convexity and this Lipschitz bound imply
\[
0\leq\kappa_a\leq\frac1{\operatorname{Rad}(\Sigma)},\qquad a=1,\ldots,n-1.
\]
Lemma~\ref{lem:metric-change} now gives
\[
\frac{n-1}{\operatorname{Rad}(\Sigma)}\leq\sqrt{H_\Sigma^2-(\tr_\Sigma k)^2}=\bar H_\Sigma=\sum_{a=1}^{n-1}\kappa_a\leq\frac{n-1}{\operatorname{Rad}(\Sigma)}.
\]
Consequently every $\kappa_a$ equals $\operatorname{Rad}(\Sigma)^{-1}$. For the Euclidean position vector $x$, we obtain
\[
d\bigl(x-\operatorname{Rad}(\Sigma)N_0(x)\bigr)=dx-\operatorname{Rad}(\Sigma)dN_0=0.
\]
Since $\partial K$ is connected, $x-\operatorname{Rad}(\Sigma)N_0(x)$ is a constant vector $c$. The map $N$ has positive degree and is onto, hence $N_0$ is onto and
\[
\partial K=\{c+\operatorname{Rad}(\Sigma)v:v\in S^{n-1}\}.
\]
By convexity, $K$ is the closed ball bounded by this sphere.	

Conversely, suppose that $f$ is constant on $\Sigma$, $k=(1+|df|_g^2)^{-1/2}\Hess_g f$, and $(M,\bar g)$ is isometric to the closed Euclidean ball of radius $\operatorname{Rad}(\Sigma)$. Lemma~\ref{lem:metric-change} gives
\[
H_\Sigma=\sqrt{(\tr_\Sigma k)^2+\bar H_\Sigma^2}=\sqrt{(\tr_\Sigma k)^2+\frac{(n-1)^2}{\operatorname{Rad}(\Sigma)^2}}\qquad\text{on }\Sigma.\qedhere
\]
\end{proof}

\section{The stable systole}\label{sec:stable-systole}

\subsection{The systolic inequality}

Suppose that the assumptions of Theorem~\ref{thm:stable-systole} hold, and let $i:\Sigma\hookrightarrow M$ be the inclusion. Write $\mathbb T^r=(\R/2\pi\Z)^r$ and $\Theta=(\Theta_1,\ldots,\Theta_{n-2})$. Let $\eta_i=(2\pi)^{-1}d\Theta_i$. Then $d\eta_i=0$ and $\int_\gamma\eta_i\in\mathbb Z$ for every closed curve $\gamma$ in $M$. Since $[\eta_i]$ are integral cohomology classes, Poincar\'e duality gives an integral homology class characterized by its pairing with closed one-forms. Writing this class as $\ell m$, we have
\begin{equation}\label{eq:systolic-pairing}
\int_\Sigma\alpha\wedge i^*\eta_1\wedge\cdots\wedge i^*\eta_{n-2}=\ell\int_m\alpha
\end{equation}
for every closed one-form $\alpha$ on $\Sigma$. Here $m$ is integral and $\ell$ is a positive integer. We can choose $\ell$ maximal since
\[
\begin{aligned}
0<|\deg(\theta_0,\Theta|_\Sigma)|&=\left|\int_\Sigma\frac{d\theta_0}{2\pi}\wedge i^*\eta_1\wedge\cdots\wedge i^*\eta_{n-2}\right|=\ell\left|\frac1{2\pi}\int_m d\theta_0\right|.
\end{aligned}
\]

\begin{lemma}\label{lem:systolic-cover-maps}
There are connected finite regular coverings $p_j:M_j\to M$ and smooth maps $\theta_j:\partial M_j\to S^1$ and $\phi_j:M_j\to\mathbb T^{n-2}$ with the following properties:
\begin{itemize}
	\item[(i)] We have
\[
\limsup_j\|d\theta_j\|_\infty\leq\frac{2\pi}{\|m\|},\qquad \|d\phi_j\|_\infty\longrightarrow0.
\]
	\item[(ii)] There is a fixed smooth degree-one map $F:\mathbb T^{n-1}\to S^{n-1}$ with $|\partial_{\theta_0}F|\leq1$.
	\item[(iii)] The maps $N_j=F(\theta_j,\phi_j|_{\partial M_j})$ have positive degree and satisfy
\[
\limsup_j\sup_{\partial M_j}\|dN_j\|_{\tr}\leq\frac{2\pi}{\|m\|}.
\]
\end{itemize}
\end{lemma}

\begin{proof}
The pairings of $m$ with integral classes in $H^1(\Sigma;\R)$ form $r\Z$ for an integer $r\geq1$. If $r>1$, then $m/r$ pairs integrally with every integral cohomology class, so $m/r$ is integral by \cite[Lemma~3.2]{BangertKatz}, contrary to maximality of $\ell$. Hence there is a smooth closed one-form $\alpha$ representing an integral class with $\int_m\alpha=1$. Fix $x_0\in\Sigma$. Since $\int_\gamma\alpha\in\Z$ for every closed curve $\gamma$, the map
\[
x\longmapsto\exp\left(2\pi\ii\int_{x_0}^x\alpha\right)
\]
is independent of the path of integration. Replacing $\theta_0$ by this smooth circle map, \eqref{eq:systolic-pairing} gives
\[
\frac1{2\pi}\int_m d\theta_0=1,\qquad \deg(\theta_0,\Theta|_\Sigma)=\int_\Sigma\frac{d\theta_0}{2\pi}\wedge i^*\eta_1\wedge\cdots\wedge i^*\eta_{n-2}=\ell>0.
\]
The classes $[d\theta_0/(2\pi)],[i^*\eta_1],\ldots,[i^*\eta_{n-2}]$ have nonzero cup product, so they are linearly independent. Since $b_1(\Sigma)=n-1$, they form a basis of $H^1(\Sigma;\R)$. Moreover,
\[
\ell\int_m i^*\eta_i=\int_\Sigma i^*\eta_i\wedge i^*\eta_1\wedge\cdots\wedge i^*\eta_{n-2}=0.
\]
Consequently, every smooth closed one-form $\alpha$ with $\int_m\alpha=1$ can be written as
\[
\alpha=\frac{d\theta_0}{2\pi}+\sum_{i=1}^{n-2}b_i\,i^*\eta_i+du,\qquad b_i\in\R.
\]

Stable norm and comass duality \cite[Section~4]{BangertKatz} give
\[
\inf\left\{\|\alpha\|_\infty:d\alpha=0,\ \int_m\alpha=1\right\}=\frac1{\|m\|},
\]
where the infimum is over smooth one-forms. We can choose real numbers $c_{ij}$ and smooth functions $u_j$ such that
\[
\left\|\frac{d\theta_0}{2\pi}+\sum_{i=1}^{n-2}c_{ij}i^*\eta_i+du_j\right\|_\infty\leq\frac1{\|m\|}+\frac1j.
\]
We next choose $b_{ij}\in\mathbb Q$ with $\sum_i|b_{ij}-c_{ij}|\,\|i^*\eta_i\|_\infty<1/j$, and set
\[
\alpha_j=\frac{d\theta_0}{2\pi}+\sum_{i=1}^{n-2}b_{ij}i^*\eta_i+du_j.
\]
Then we have $\int_m\alpha_j=1$ and
\[
\frac1{\|m\|}\leq\|\alpha_j\|_\infty\leq\frac1{\|m\|}+\frac2j.
\]
Choose positive integers $L_j\geq j$ such that $L_jb_{ij}\in\mathbb Z$ for $i=1,\ldots,n-2$, and let $M_j$ be a connected component of
\[
\{(x,z)\in M\times\mathbb T^{n-2}:\Theta(x)=L_jz\}.
\]
Let $p_j:M_j\to M$ and $\phi_j=(\phi_{1,j},\ldots,\phi_{n-2,j}):M_j\to\mathbb T^{n-2}$ be the projections. The covering $p_j$ corresponds to the subgroup $\Theta_*^{-1}(L_j\Z^{n-2})$, which is normal and has finite index in $\pi_1(M)$. Thus $p_j$ is a connected finite regular covering. We equip $M_j$ with the pullbacks of $g$, $k$, and the spin structure. The maps satisfy
\[
L_j\phi_j=\Theta\circ p_j,\qquad d\phi_{i,j}=\frac{2\pi}{L_j}p_j^*\eta_i.
\]
On $\partial M_j$, define
\[
\theta_j=\theta_0\circ p_j+\sum_i(L_jb_{ij})\phi_{i,j}+2\pi u_j\circ p_j.
\]
Since $L_jb_{ij}\in\Z$, this defines a smooth circle map. Then
\[
d\theta_j=p_j^*d\theta_0+\sum_i(L_jb_{ij})d\phi_{i,j}+2\pi p_j^*du_j=2\pi p_j^*\alpha_j.
\]
Since $p_j$ is a surjective local isometry and $d\phi_j=L_j^{-1}p_j^*d\Theta$, we obtain
\[
\begin{aligned}
\|d\theta_j\|_\infty&=2\pi\|p_j^*\alpha_j\|_\infty=2\pi\|\alpha_j\|_\infty\leq\frac{2\pi}{\|m\|}+\frac{4\pi}{j},\\
\|d\phi_j\|_\infty&=\frac1{L_j}\|p_j^*d\Theta\|_\infty=\frac1{L_j}\|d\Theta\|_\infty\longrightarrow0.
\end{aligned}
\]
This proves (i).

The degree on $\partial M_j$ is the sum of the degrees of the restrictions to its connected components, with their induced boundary orientations. By \eqref{eq:systolic-pairing},
\[
\begin{aligned}
\deg(\theta_j,\phi_j|_{\partial M_j})&=\frac1{(2\pi)^{n-1}}\int_{\partial M_j}d\theta_j\wedge d\phi_{1,j}\wedge\cdots\wedge d\phi_{n-2,j}\\
&=\frac{\deg p_j}{L_j^{n-2}}\int_\Sigma\alpha_j\wedge i^*\eta_1\wedge\cdots\wedge i^*\eta_{n-2}\\
&=\frac{\ell\deg p_j}{L_j^{n-2}}\\
&>0.
\end{aligned}
\]

Iterating on the construction in the proof of \cite[Lemma~3.8]{BBCH} gives a smooth degree-one map $F:\mathbb T^{n-1}\to S^{n-1}$. In particular, 
\[
|\partial_{\theta_0}F|\leq|\partial_{\theta_0}(\cos\theta_0,\sin\theta_0)|=1.
\]
Since $F$ is smooth on the compact torus, $|\partial_{\theta_i}F|\leq C$ for $1\leq i\leq n-2$. Let $N_j=F(\theta_j,\phi_j|_{\partial M_j})$. Then $\deg N_j>0$, and
\[
\begin{aligned}
\|dN_j\|_{\tr}&=\left\|\partial_{\theta_0}F\otimes d\theta_j+\sum_{i=1}^{n-2}\partial_{\theta_i}F\otimes d\phi_{i,j}\right\|_{\tr}\\
&\leq|\partial_{\theta_0}F|\,|d\theta_j|+\sum_{i=1}^{n-2}|\partial_{\theta_i}F|\,|d\phi_{i,j}|\\
&\leq2\pi|p_j^*\alpha_j|+\frac{2\pi C}{L_j}\sum_{i=1}^{n-2}|p_j^*\eta_i|.
\end{aligned}
\]
This gives
\[
\begin{aligned}
\sup_{\partial M_j}\|dN_j\|_{\tr}&\leq2\pi\|\alpha_j\|_\infty+\frac{2\pi C}{L_j}\sum_{i=1}^{n-2}\|\eta_i\|_\infty
\leq\frac{2\pi}{\|m\|}+\frac{4\pi}{j}+\frac{2\pi C}{L_j}\sum_{i=1}^{n-2}\|\eta_i\|_\infty.
\end{aligned}
\]
Letting $j\to\infty$ proves (iii).
\end{proof}

\begin{proof}[Proof of the inequality in Theorem~\ref{thm:stable-systole}]
Suppose that
\[
\min_\Sigma\left(H_\Sigma-\sqrt{(\tr_\Sigma k)^2+\frac{4\pi^2}{\|m\|^2}}\right)>0.
\]
By Lemma~\ref{lem:systolic-cover-maps}(iii), there are coverings $p_j:M_j\to M$, positive-degree maps $N_j:\partial M_j\to S^{n-1}$, and $\varepsilon_j>0$ with $\varepsilon_j\to0$ such that $\sup_{\partial M_j}\|dN_j\|_{\tr}\leq2\pi/\|m\|+\varepsilon_j$. For sufficiently large $j$, we have
\[
\sqrt{(\tr_\Sigma k)^2+\|dN_j\|_{\tr}^2}\leq\sqrt{(\tr_\Sigma k)^2+\frac{4\pi^2}{\|m\|^2}}+\varepsilon_j<H_\Sigma.
\]
Theorem~\ref{thm:boundary-rigidity} gives
\[
H_\Sigma=\sqrt{(\tr_\Sigma k)^2+\|dN_j\|_{\tr}^2}<H_\Sigma
\]
on $\partial M_j$, a contradiction. Since $m$ is a nonzero integral class, $\operatorname{stsys}_1(\Sigma)\leq\|m\|$, and hence
\begin{equation}\label{eq:systolic-selected-bound}
\begin{aligned}
\min_\Sigma\left(H_\Sigma-\sqrt{(\tr_\Sigma k)^2+\frac{4\pi^2}{\operatorname{stsys}_1(\Sigma)^2}}\right)&\leq\min_\Sigma\left(H_\Sigma-\sqrt{(\tr_\Sigma k)^2+\frac{4\pi^2}{\|m\|^2}}\right)\leq0.
\end{aligned}
\end{equation}
\end{proof}

\subsection{The equality case}

\begin{lemma}\label{lem:systolic-current}
Suppose that $n$ is odd and
\[
\min_\Sigma\left(H_\Sigma-\sqrt{(\tr_\Sigma k)^2+\frac{4\pi^2}{\|m\|^2}}\right)=0.
\]
Then there is a vector field $V$ on $M$ with the following properties:
\begin{itemize}
	\item[(i)] We have
\[
D_XV=-\sqrt{1+|V|_g^2}\,k(X),\qquad V^\flat|_{T\Sigma}=0.
\]
	\item[(ii)] The metric $\bar g=g+V^\flat\otimes V^\flat$ is Ricci flat and induces the same metric on $\Sigma$ as $g$.
	\item[(iii)] The second fundamental form $\bar h$ and mean curvature $\bar H_\Sigma$ of $\Sigma$ in $(M,\bar g)$ satisfy
\[
\bar h\geq0,\qquad \bar H_\Sigma=\frac{2\pi}{\|m\|}.
\]
\end{itemize}
\end{lemma}

\begin{proof}
Let $\pi:\widetilde M\to M$ be the universal cover. Since $\widetilde M$ is simply connected and each $M_j$ is connected, lifting $\pi$ through $p_j$ gives a covering map $q_j:\widetilde M\to M_j$ with $p_j\circ q_j=\pi$. We can find $\varepsilon_j>0$ with $\varepsilon_j\to0$ and $\sup_{\partial M_j}\|dN_j\|_{\tr}\leq2\pi/\|m\|+\varepsilon_j$. Then
\[
\begin{aligned}
\sqrt{(\tr_\Sigma k)^2+\|dN_j\|_{\tr}^2}&\leq\sqrt{(\tr_\Sigma k)^2+\frac{4\pi^2}{\|m\|^2}}+\varepsilon_j
\leq H_\Sigma+\varepsilon_j.
\end{aligned}
\]
By Proposition~\ref{prop:parallel-section}, the maps $q_j$ can be chosen so that $N_j\circ q_j$ converges locally uniformly to a smooth map $N:\partial\widetilde M\to S^{n-1}$, and there is a smooth nowhere-zero section $s$ satisfying $\nabla^{\mathcal E,k}s=0$ on $\widetilde M$ and $\chi s=s$ on $\partial\widetilde M$.

Since $\widetilde M$ is connected and $\|d\phi_j\|_\infty\to0$ by Lemma~\ref{lem:systolic-cover-maps}(i), after passing to a subsequence, $\phi_j\circ q_j$ converges locally uniformly to a constant $\phi_*\in\mathbb T^{n-2}$. By the same lemma, after passing to a further subsequence, $\theta_j\circ q_j$ converges locally uniformly on $\partial\widetilde M$ to a Lipschitz map $\theta$ with $|d\theta|\leq2\pi/\|m\|$ almost everywhere. Hence $N=F(\theta,\phi_*)$, and
\[
dN=\partial_{\theta_0}F(\theta,\phi_*)\otimes d\theta
\]
almost everywhere. Proposition~\ref{prop:parallel-section} and Lemma~\ref{lem:systolic-cover-maps}(ii) give
\[
\begin{aligned}
\frac{2\pi}{\|m\|}&\leq\sqrt{H_\Sigma^2-(\tr_\Sigma k)^2}\\
&=\|dN\|_{\tr}\\
&=|\partial_{\theta_0}F(\theta,\phi_*)|\,|d\theta|\\
&\leq|d\theta|\\
&\leq\frac{2\pi}{\|m\|}.
\end{aligned}
\]
almost everywhere. Since $N$ is smooth, $H_\Sigma^2-(\tr_\Sigma k)^2=4\pi^2/\|m\|^2$ and $dN$ has rank one everywhere on $\partial\widetilde M$.

We next define the vector field $V$ by \eqref{eq:current-definition}. Lemma~\ref{lem:current} gives $V=g(V,\nu)\nu$ on $\partial\widetilde M$ and
\[
H_\Sigma g(V,\nu)=\sum_a g(D_{e_a}V,e_a)=-|s|^2\tr_\Sigma k.
\]
Consequently,
\[
|s|^4-|V|_g^2=\frac{|s|^4}{H_\Sigma^2}\bigl(H_\Sigma^2-(\tr_\Sigma k)^2\bigr)=\frac{4\pi^2|s|^4}{\|m\|^2H_\Sigma^2}>0.
\]
on $\partial\widetilde M$. By Lemma~\ref{lem:current}, $|s|^4-|V|_g^2$ is constant on $\widetilde M$. After rescaling, we may assume that $|s|^4-|V|_g^2=1$. Thus, we have
\[
|s|^2=\frac{\|m\|H_\Sigma}{2\pi},\qquad V=-\frac{\|m\|\tr_\Sigma k}{2\pi}\nu.
\]
on $\partial\widetilde M$. Let $\tau$ be a deck transformation of $\pi$. Since $\pi\circ\tau=\pi$, we have $\tau^*V=V$ and $\tau^*|s|^2=|s|^2$ on $\partial\widetilde M$. Lemma~\ref{lem:current} and ODE uniqueness then give $\tau^*V=V$ and $\tau^*|s|^2=|s|^2$ on $\widetilde M$. Any two lifts of $p\in M$ are related by a deck transformation, so $d\pi_x(V_x)$ and $|s(x)|^2$ are independent of the choice of $x\in\pi^{-1}(p)$. They define a vector field $V$ and a function $|s|^2$ on $M$, where
\[
D_XV=-|s|^2k(X)=-\sqrt{1+|V|_g^2}\,k(X),\qquad V^\flat|_{T\Sigma}=0.
\]
This proves (i).

Proposition~\ref{prop:comparison-metric}, applied on $\widetilde M$, gives (ii) and
\[
\bar H_\Sigma=\frac{H_\Sigma}{|s|^2}=\frac{2\pi}{\|m\|}>0.
\]
It also gives $\ker dN\subset\ker\bar h$. Since $dN$ has rank one, $\bar h$ has rank at most one. Its positive trace implies $\bar h\geq0$, proving (iii).
\end{proof}

\begin{lemma}\label{lem:systolic-disk-product}
Under the assumptions of Theorem~\ref{thm:stable-systole}, let $\bar g$ be a Ricci-flat metric on $M$ with $\bar h\geq0$ and constant $\bar H_\Sigma>0$. Then $i^*:H^1(M;\R)\to H^1(\Sigma;\R)$ is injective, and the universal cover of $(M,\bar g)$ is isometric to
\[
B^2(1/\bar H_\Sigma)\times\R^{n-2}
\]
with its Euclidean product metric.
\end{lemma}

\begin{proof}
Let $\eta_i=(2\pi)^{-1}d\Theta_i$, $i=1,\ldots,n-2$. Since $\int_Md^*\eta_i=-\int_\Sigma\eta_i(\bar\nu)$, the Neumann problems
\[
\Delta_{\bar g}u_i=d^*\eta_i\quad\text{on }M,\qquad \partial_{\bar\nu}u_i=-\eta_i(\bar\nu)\quad\text{on }\Sigma
\]
have smooth solutions. Then $\omega_i=\eta_i+du_i$ represents $[\eta_i]$ and satisfies $d\omega_i=d^*\omega_i=0$ and $\omega_i(\bar\nu)=0$. Writing $\bar D$ for the Levi--Civita connection of $\bar g$, the Weitzenb\"ock formula and $\Ric_{\bar g}=0$ give
\[
\begin{aligned}
0&=\int_M\langle(dd^*+d^*d)\omega_i,\omega_i\rangle\\
&=\int_M|\bar D\omega_i|^2-\int_\Sigma(\bar D_{\bar\nu}\omega_i)(\omega_i^\sharp)\\
&=\int_M|\bar D\omega_i|^2-\int_\Sigma(\bar D_{\omega_i^\sharp}\omega_i)(\bar\nu)\\
&=\int_M|\bar D\omega_i|^2+\int_\Sigma\bar h(\omega_i^\sharp,\omega_i^\sharp).
\end{aligned}
\]
Thus each $\omega_i$ is parallel, and $\omega_i^\sharp$ is tangent to $\Sigma$ and lies in $\ker\bar h$. Since $\deg(\theta_0,\Theta|_\Sigma)\neq0$, $[\eta_1],\ldots,[\eta_{n-2}]$ are linearly independent. A linear combination of the parallel $\omega_i$ that vanishes at one point vanishes on $M$, so the $\omega_i$ are pointwise linearly independent. After a constant linear change of their duals, we obtain orthonormal parallel vector fields $E_1,\ldots,E_{n-2}$ tangent to $\Sigma$ and lying in $\ker\bar h$. Since $\bar R(X,Y)E_a=0$ for every $a$, all components involving an $E_a$ vanish. For an orthonormal basis $e_1,e_2$ of their orthogonal complement, we have
\[
\begin{aligned}
0=\Ric_{\bar g}(e_1,e_1)&=\bar g(\bar R(e_2,e_1)e_1,e_2)+\sum_{a=1}^{n-2}\bar g(\bar R(E_a,e_1)e_1,E_a)
=\bar g(\bar R(e_2,e_1)e_1,e_2).
\end{aligned}
\]
Thus $\bar g$ is flat.

Let $\pi:\widetilde M\to M$ be the universal cover. Lifting the metric $\bar g$ to  $\widetilde M$ gives a local isometry $G:\widetilde M\to\R^n$. The lifted metric is complete, and $\bar h\geq0$ makes the boundary locally convex. For distinct $p,q\in\widetilde M$, let $\gamma:[0,L]\to\widetilde M$ be a minimizing curve joining them \cite[Theorem~2.5.23]{BBI}, parametrized by arc length. Each sufficiently short subarc lies in a neighbourhood on which $G$ is injective with convex image and therefore maps to a straight line segment by minimality of $\gamma$. Hence $(G\circ\gamma)''=0$ and $|(G\circ\gamma)'|=1$, which gives
\[
\begin{aligned}
|G(q)-G(p)|&=\left|\int_0^L(G\circ\gamma)'(t)\,dt\right|\\
&=L|(G\circ\gamma)'(0)|\\
&=\operatorname{Length}(G\circ\gamma)\\
&=d_{\widetilde M}(p,q).
\end{aligned}
\]
Thus $G$ is injective, its image is convex, and is also closed by completeness. We then identify $\widetilde M$ with this closed convex domain $K\subset\R^n$.

The parallel fields $E_1,\ldots,E_{n-2}$ become constant vectors tangent to $\partial K$. Every supporting half-space of $K$ has normal orthogonal to their span and is therefore invariant under translations in these directions. Since $K$ is the intersection of its supporting half-spaces, after an orthogonal change of coordinates we have
\[
K=K_0\times\R^{n-2}
\]
for a closed convex planar domain $K_0$. Its boundary is smooth and has constant curvature $\bar H_\Sigma>0$. Each component is closed in $\R^2$ and has constant positive curvature, hence it is a circle of radius $1/\bar H_\Sigma$. Convexity makes $K_0$ the disk bounded by this circle. Thus $K_0=B^2(1/\bar H_\Sigma)$, up to translation.

Finally, let $\alpha$ be a closed one-form on $M$ with $i^*[\alpha]=0$. Fix $x_0\in\Sigma$ and a loop $\gamma$ in $M$ based at $x_0$. A lift $\widetilde\gamma$ has both endpoints in $\pi^{-1}(x_0)\subset\partial\widetilde M$. Since $\partial\widetilde M=S^1(1/\bar H_\Sigma)\times\R^{n-2}$ is connected, there is a path $\widetilde\sigma$ in $\partial\widetilde M$ joining these endpoints. The form $\pi^*\alpha$ is exact on the simply connected manifold $\widetilde M$, and $\pi\circ\widetilde\sigma$ is a loop in $\Sigma$. Hence
\[
\int_\gamma\alpha=\int_{\widetilde\gamma}\pi^*\alpha=\int_{\widetilde\sigma}\pi^*\alpha=\int_{\pi\circ\widetilde\sigma}i^*\alpha=0.
\]
Thus $\alpha$ is exact, proving that $i^*:H^1(M;\R)\to H^1(\Sigma;\R)$ is injective.
\end{proof}

\begin{proof}[Proof of the equality statement in Theorem~\ref{thm:stable-systole}]
Suppose that equality holds in Theorem~\ref{thm:stable-systole}. By \eqref{eq:systolic-selected-bound},
\[
\min_\Sigma\left(H_\Sigma-\sqrt{(\tr_\Sigma k)^2+\frac{4\pi^2}{\|m\|^2}}\right)=0.
\]
Let $p\in\Sigma$ be a point where this minimum is attained. Then
\[
\begin{aligned}
\sqrt{(\tr_\Sigma k(p))^2+\frac{4\pi^2}{\operatorname{stsys}_1(\Sigma)^2}}&\leq H_\Sigma(p)\\
&=\sqrt{(\tr_\Sigma k(p))^2+\frac{4\pi^2}{\|m\|^2}}\\
&\leq\sqrt{(\tr_\Sigma k(p))^2+\frac{4\pi^2}{\operatorname{stsys}_1(\Sigma)^2}}.
\end{aligned}
\]
Thus $\|m\|=\operatorname{stsys}_1(\Sigma)$.

If $n$ is odd, define $V$ and $\bar g=g+V^\flat\otimes V^\flat$ as in Lemma~\ref{lem:systolic-current}. The one form $V^\flat$ is closed by Lemma~\ref{lem:current} and vanishes on $\Sigma$. Applying Lemma~\ref{lem:systolic-disk-product}, we can find a smooth function $f$ with $V^\flat=-df$, constant on  $\Sigma$, and such that the universal cover of $(M,\bar g)$ is isometric to $B^2((2\pi)^{-1}\operatorname{stsys}_1(\Sigma))\times\R^{n-2}$. Moreover,
\[
\begin{aligned}
\Hess_g f(X,Y)&=-(D_XV^\flat)(Y)=\sqrt{1+|V|_g^2}\,k(X,Y)
=\sqrt{1+|df|_g^2}\,k(X,Y).
\end{aligned}
\]

If $n$ is even, we consider the initial data set $(M\times S^1(L),g+dt^2,k\oplus0)$ with $L\geq\operatorname{stsys}_1(\Sigma)$, where $t$ denotes the arc-length parameter on $S^1(L)$. Note that $b_1(\Sigma\times S^1(L))=n$ and it satisfies the assumptions of Theorem~\ref{thm:stable-systole} in odd dimension. Moreover, the map $(x,t)\mapsto(\theta_0(x),\Theta(x),e^{\ii t/L})$ from $\Sigma\times S^1(L)$ to $S^1\times\mathbb T^{n-1}$ has nonzero degree, and its last $n-1$ components extend to $M\times S^1(L)$. Under the identification $H_1(\Sigma\times S^1(L);\R)=H_1(\Sigma;\R)\oplus\R$, the integral classes are $(\alpha,j)$ with $\alpha$ integral and $j\in\Z$. The projections to the two factors are $1$-Lipschitz, so
\[
\|(\alpha,j)\|\geq\max\{\|\alpha\|,2\pi L|j|\}\geq\operatorname{stsys}_1(\Sigma)
\]
if $(\alpha,j)\neq0$. The inclusion $x\mapsto(x,0)$ is also $1$-Lipschitz, hence
\[
\operatorname{stsys}_1(\Sigma)\leq\operatorname{stsys}_1(\Sigma\times S^1(L))\leq\operatorname{stsys}_1(\Sigma).
\]
Since the product also has the same boundary mean curvature and tangential trace of $k$, it satisfies equality in Theorem~\ref{thm:stable-systole} in odd dimension. Consequently, we can find a smooth function $f$, constant on $\Sigma\times S^1(L)$, such that
\[
\Hess_{g+dt^2}f=\sqrt{1+|df|_{g+dt^2}^2}\,(k\oplus0).
\]
In particular,
\[
\partial_t^2f=\Hess_{g+dt^2}f(\partial_t,\partial_t)=\sqrt{1+|df|_{g+dt^2}^2}\,(k\oplus0)(\partial_t,\partial_t)=0.
\]
Pulling $f$ back to $M\times\R$, we obtain
\[
0=f(x,t+2\pi L)-f(x,t)=\int_t^{t+2\pi L}\partial_tf(x,s)\,ds=2\pi L\,\partial_tf(x,t).
\]
Thus $f$ is independent of $t$, and restriction to $M$ gives $k=(1+|df|_g^2)^{-1/2}\Hess_g f$. In particular, $g+df\otimes df$ is flat, its boundary second fundamental form is nonnegative, and its boundary mean curvature is $2\pi/\operatorname{stsys}_1(\Sigma)$. Lemma~\ref{lem:systolic-disk-product} then gives the required universal cover of $(M,g+df\otimes df)$.

Conversely, let $f\in C^\infty(M)$ be constant on $\Sigma$, and let $\bar g=g+df\otimes df$. Suppose that $k=(1+|df|_g^2)^{-1/2}\Hess_g f$ and the universal cover of $(M,\bar g)$ is isometric to $B^2((2\pi)^{-1}\operatorname{stsys}_1(\Sigma))\times\R^{n-2}$ with the Euclidean product metric. Then $\bar g|_\Sigma=g|_\Sigma$ and $\bar H_\Sigma=2\pi/\operatorname{stsys}_1(\Sigma)$. By Lemma~\ref{lem:metric-change},
\[
H_\Sigma=\sqrt{(\tr_\Sigma k)^2+\bar H_\Sigma^2}=\sqrt{(\tr_\Sigma k)^2+\frac{4\pi^2}{\operatorname{stsys}_1(\Sigma)^2}}.\qedhere
\]
\end{proof}



\begin{thebibliography}{99}
\bibitem{Ambrozio}
L. C. Ambrozio, \emph{Rigidity of area-minimizing free boundary surfaces in mean convex three-manifolds}, J. Geom. Anal. \textbf{25} (2015), 1001--1017.

\bibitem{AmmannLockman}
B. Ammann and S. Lockman, \emph{Rigidity for spin fill-ins with scalar curvature bounded from below}, arXiv:2608.06935 (2026).

\bibitem{BarRadius}
C. B\"ar, \emph{Dirac eigenvalues and the hyperspherical radius}, J. Eur. Math. Soc. (2026), published online.

\bibitem{BangertKatz}
V. Bangert and M. Katz, \emph{Stable systolic inequalities and cohomology products}, Comm. Pure Appl. Math. \textbf{56} (2003), 979--997.

\bibitem{BBCH}
C. B\"ar, S. Brendle, T.-K. A. Chow, and B. Hanke, \emph{Rigidity results for initial data sets satisfying the dominant energy condition}, J. reine angew. Math. \textbf{836} (2026), 27--69.

\bibitem{BGM}
C. B\"ar, P. Gauduchon, and A. Moroianu, \emph{Generalized cylinders in semi-Riemannian and spin geometry}, Math. Z. \textbf{249} (2005), 545--580.

\bibitem{BrendleHung}
S. Brendle and P.-K. Hung, \emph{Systolic inequalities and the Horowitz--Myers conjecture}, Acta Math., to appear; arXiv:2406.04283, version~4 (2026).

\bibitem{BMN}
S. Brendle, F. C. Marques, and A. Neves, \emph{Deformations of the hemisphere that increase scalar curvature}, Invent. Math. \textbf{185} (2011), 175--197.

\bibitem{BrendleTsiamisWang}
S. Brendle, R. Tsiamis, and Y. Wang, \emph{On fill-ins with scalar curvature bounded from below and an inequality of Hijazi--Montiel--Rold\'an}, arXiv:2510.17780 (2025), version~2.

\bibitem{BBI}
D. Burago, Y. Burago, and S. Ivanov, \emph{A Course in Metric Geometry}, Graduate Studies in Mathematics, vol.~33, American Mathematical Society, Providence, RI, 2001.

\bibitem{CHZ}
S. Cecchini, S. Hirsch, and R. Zeidler, \emph{Rigidity of spin fill-ins with non-negative scalar curvature}, arXiv:2404.17533 (2024), version~3 (2026).

\bibitem{GromovLectures}
M. Gromov, \emph{Four lectures on scalar curvature}, arXiv:1908.10612 (2019), version~6 (2021).

\bibitem{HMR}
O. Hijazi, S. Montiel, and A. Rold\'an, \emph{Dirac operators on hypersurfaces of manifolds with negative scalar curvature}, Ann. Global Anal. Geom. \textbf{23} (2003), 247--264.

\bibitem{HuLiuShi}
Y. Hu, P. Liu, and Y. Shi, \emph{Rigidity of 3D spherical caps via $\mu$-bubbles}, Pacific J. Math. \textbf{323} (2023), 89--114.

\bibitem{KoYao}
D. Ko and X. Yao, \emph{Capillary minimal slicing and scalar curvature rigidity}, arXiv:2602.21071 (2026), version~2.

\bibitem{LeeTam}
M.-C. Lee and L.-F. Tam, \emph{Rigidity of Lipschitz map using harmonic map heat flow}, Amer. J. Math. \textbf{148} (2026), 985--1002.

\bibitem{Llarull}
M. Llarull, \emph{Sharp estimates and the Dirac operator}, Math. Ann. \textbf{310} (1998), 55--71.

\bibitem{Longa}
E. Longa, \emph{Sharp systolic inequalities for $3$-manifolds with boundary}, J. Geom. Anal. \textbf{31} (2021), 7741--7749.

\bibitem{LottBoundary}
J. Lott, \emph{Index theory for scalar curvature on manifolds with boundary}, Proc. Amer. Math. Soc. \textbf{149} (2021), 4451--4459.

\bibitem{RaulotDEC}
S. Raulot, \emph{Nonexistence of DEC spin fill-ins}, C. R. Math. Acad. Sci. Paris \textbf{360} (2022), 1049--1054.

\bibitem{RaulotFillin}
S. Raulot, \emph{Dirac eigenvalues and spacetime mean curvature estimates for DEC spin fill-ins}, arXiv:2609.08720 (2026).

\bibitem{Wang}
Y. Wang, \emph{Rigidity results involving stabilized scalar curvature}, arXiv:2501.06951 (2025).

\bibitem{WangXie}
J. Wang and Z. Xie, \emph{On Gromov's flat corner domination conjecture and Stoker's conjecture}, arXiv:2203.09511 (2022), version~2 (2023).

\bibitem{YaoComparison}
X. Yao, \emph{A note on scalar curvature comparison rigidity for compact domains}, arXiv:2410.21238 (2024).
\end{thebibliography}
\end{document}